\documentclass{article}

\usepackage{amsfonts}

\usepackage{amssymb}

\newcommand{\marcio}[1]{\stackrel{#1}{_{\longmapsto}}}
\newcommand{\rank}[0]{\mathtt{rank}}

\title{Flow: the Axiom of Choice is independent from the Partition Principle in ZFU}

\author{Adonai S. Sant'Anna \and Renato Brodzinski \and Marcio P. P. de Fran\c ca \and Ot\'avio Bueno}

\date{adonaisantanna@gmail.com, renatobrodzinski@gmail.com, marciopalmares@gmail.com, otaviobueno@me.com}

\begin{document}

\newtheorem{definicao}{Definition}
\newtheorem{teorema}{Theorem}
\newtheorem{lema}{Lemma}
\newtheorem{corolario}{Corollary}
\newtheorem{proposicao}{Proposition}
\newtheorem{axioma}{Axiom}
\newtheorem{observacao}{Observation}




\maketitle

\begin{abstract}

We introduce a formal theory called Flow where the intended interpretation of its terms is that of {\em function\/}. We prove ZF, ZFC  and ZFU (ZF with atoms) can be immersed within Flow as natural consequences from our framework. Our first important application is the introduction of a model of ZFU where the Partition Principle holds but the Axiom of Choice fails, if Flow is consistent. So, our framework allows us to address the oldest open problem in set theory: if the Partition Principle entails the Axiom of Choice.\\

\begin{center}
{\sc\large This is a fully revised version of a previous preprint about Flow and its applications.}
\end{center}
\end{abstract}

\tableofcontents

\section{Introduction}

It is rather difficult to  determine when functions were born, since mathematics itself significantly evolves along history. Specially nowadays we find different formal concepts associated to the label {\em function\/}. But an educated guess could point towards Sharaf al-D\={\i}n al-T\={u}s\={\i} who, in the 12th century, not only introduced a `dynamical' concept which could be interpreted as some notion of function, but also studied how to determine the maxima of such functions \cite{Hogenduk-89}.

All usual mathematical approaches for physical theories are based on either differential equations or systems of differential equations whose solutions (when they exist) are either functions or classes of functions (see, e.g., \cite{Sobolev-11}). In pure mathematics the situation is no different. Continuous functions, linear transformations, homomorphisms, and homeomorphisms, for example, play a fundamental role in topology, linear algebra, group theory, and differential geometry, respectively. Category theory \cite{MacLane-94} emphasizes such a role in a very clear, elegant, and comprehensive way. As remarked by Marquis \cite{Marquis-18}, category theory allows us to distinguish between canonical and noncanonical maps, in a way which is not usually achieved within purely extensional set theories, like ZFC. And canonical maps ``constitute the highway system of mathematical concepts''. Concepts of symmetry are essential in pure and applied mathematics, and they are stated by means of group transformations \cite{Sternberg-95}.

Some authors suggest that functions are supposed to play a strategic role into the foundations of mathematics \cite{vonNeumann-25} and even mathematics teaching \cite{Klein-16}, rather than sets. The irony of such perspective starts in the historical roots of set theories. Georg Cantor's seminal works on sets were strongly motivated by Bernard Bolzano's manuscripts on infinite multitudes called {\em Menge\/} \cite{Trlifajova-18}. Those collections were supposed to be conceived in a way such that the arrangement of their components is unimportant. However, Bolzano insisted on an Euclidian view that the whole should be greater than a part, while Cantor proposed a quite different approach: to compare infinite quantities we should consider a {\em one-to-one correspondence\/} between collections. Cantor's concept of collection (in his famous {\em Mengenlehre\/}) was strongly committed to the idea of function. Subsequent formalizations of Cantor's ``theory'' were developed in a way such that all strategic terms were associated to an intended interpretation of collection. The result of that effort is a strange phenomenon from the point of view of theories of definition: Padoa's principle allows us to show domains of functions are definable from functions themselves \cite{daCosta-01} \cite{daCosta-02} \cite{Sant'Anna-14}. Thus, even in extensional set theories functions seem to play a more fundamental role than sets.

Nevertheless, there is another crucial role that functions (in a broader intuitive sense) play in the study of set theories and cannot be ignored, namely, those questions regarding the use of $\in$-automorphisms \cite{Jech-03} and non-trivial elementary embeddings of a universe into itself \cite{Kanamori-03}, among other examples. Some questions regarding set theories cannot be answered by means of the use of their axioms alone. That is where model theory provides metamathematical tools to answer to some of those questions \cite{Suzuki-99}. Non-trivial elementary embeddings and $\in$-automorphisms are examples of `functions' which are not formalized into standard set theories but which are very useful to evaluate questions regarding independence of formulas like the Axiom of Choice, the Continuum Hypothesis and the existence of inaccessible cardinals. Within this context, the formal theory we introduce here is supposed to provide a framework where at least some of those `very large functions' can be precisely defined.

This paper was initially motivated by von Neumann's original ideas \cite{vonNeumann-25} and a variation of them \cite{Sant'Anna-14}, and related papers as well (\cite{daCosta-01} \cite{daCosta-02}). In \cite{Sant'Anna-14} it was provided a reformulation of von Neumann's `functions' theory (termed $\mathcal N$ theory). Nevertheless, we develop here a whole new approach called {\em Flow\/}.

String diagrams \cite{Dixon-13} (whose edges and vertices can eventually be interpreted as morphisms in a monoidal category) introduce some ideas which seem to intersect with our own. A remarkable feature of string diagrams is that edges need not be connected to vertices at both ends. More than that, unconnected ends can be interpreted as inputs and outputs of a string diagram (with important applications in computer science). But within our framework {\em no\/} function has any domain whatsoever. Besides, the usual way to cope with string diagrams is by means of a discrete and finitary framework, while in Flow (our framework) we do not impose such restrictions.

A more radical proposal where functions play a fundamental role is the Theory of Autocategories \cite{Guitart-14}. Autocategories are developed with the aid of autographs, where arrows (which work as morphisms) are drawn between arrows with no need of objects. Thus, once more we see a powerful idea concerning functions is naturally emerging in different places nowadays.

When Ernst Zermelo introduced AC, his motivation was the Partition Principle \cite{Banaschewski-90}, which is usually stated by means of functions. Indeed, AC entails PP. But since 1904 (when Zermelo introduced AC) it is unknown whether PP implies AC. For a review about the subject see, e.g., \cite{Banaschewski-90} \cite{Higasikawa-95} \cite{Howard-16} \cite{Pelc-78}. Our answer is that AC is independent from PP, at least within ZFU if Flow is consistent. Our point is that ZF, ZFC, all their variants, and almost all their respective models (with a few exceptions like those in \cite{Abian-78} and \cite{Devlin-93}) are somehow committed to a methodological and epistemological character which forces us to see sets as collections of some sort. Within Flow that does not happen, since our framework drives us to see sets as special cases of terms whose intended interpretation correspond to an intuitive notion of function. All sets in Flow are restrictions of $\underline{\mathfrak{1}}$, a function which works as some sort of `universal' identity in the sense that for any $x$ we have $\underline{\mathfrak{1}}(x) = x$. Nevertheless, functions who are no restrictions of $\underline{\mathfrak{1}}$ play an important role with consequences over sets. An analogous situation takes place with the well known Reflexion Principle in model theory \cite{Steel-07}: some features of the von Neumann Universe motivate mathematicians to look for new axioms which grant the existence of strongly inaccessible cardinals. Flow, however, provides a new way for coping with the metamathematics of ZF and its variants.

As a final remark before we start, it is worth to observe we freely employ the term `Flow' in this paper. Within a more general viewpoint, `Flow' refers to any first order theory {\sc (i)} whose objects are supposed to be interpreted as functions; {\sc (ii)} which grants the existence of functions $\underline{\mathfrak{0}}$ and $\underline{\mathfrak{1}}$, in the sense that $\forall x(\underline{\mathfrak{1}}(x) = x \wedge \underline{\mathfrak{0}}(x) = \underline{\mathfrak{0}})$ and $\underline{\mathfrak{0}}\neq\underline{\mathfrak{1}}$, as it follows in the next Section; and {\sc (iii)} which allows a nontrivial composition among any functions at all, whether such a composition is associative or not. In this first paper about Flow, however, we explore one specific formal framework who is able to replicate ZF, ZFC, and ZFU. Our goal here is to suggest another way to construct models of set theories.

\section{Flow theory}\label{secao2}

Our first Flow theory is a first-order theory $\mbox{\boldmath{$\mathfrak{F}$}}$ with identity \cite{Mendelson-97}, where $x=y$ reads `$x$ is equal to $y$', and $\neg(x = y)$ is abbreviated as $x\neq y$. From now on we use Flow and $\mbox{\boldmath{$\mathfrak{F}$}}$ as synonyms. Flow has one functional letter $f_1^2$ (termed {\em evaluation\/}), where $f_1^2(f,x)$ is a term, if $f$ and $x$ are terms. If $y = f_1^2(f,x)$, we abbreviate this by $f(x) = y$, and read `$y$ is the {\em image\/} of $x$ by $f$'. All terms are called {\em functions\/}. Such a terminology seems to be adequate under the light of our intended interpretation: functions are supposed to be terms which `transform' terms into other terms. Since we are assuming identity, our functions cannot play the role of non-trivial relations. Suppose, for example, $f(x) = y$ and $f(x) = y'$, which are abbreviations for $f_1^2(f,x) = y$ and $f_1^2(f,x) = y'$, respectively. From transitivity of identity, we have $y = y'$. So, for any $f$ and any $x$, there is one single $y$ such that $f(x) = y$.

Lowercase Latin and Greek letters denote functions, with the sole exception of two specific functions to appear in the next pages, namely, $\underline{\mathfrak{0}}$ and $\underline{\mathfrak{1}}$. Uppercase Latin letters are used to denote formulas or predicates (which are eventually defined). Any explicit definition in Flow is an abbreviative one, in the sense that for a given formula $F$, the {\em definiendum\/} is a metalinguistic abbreviation for the {\em definiens\/} given by $F$. Eventually we use bounded quantifiers. If $P$ is a predicate defined by a formula $F$, we abbreviate $\forall x(P(x)\Rightarrow G(x))$ and $\exists x(P(x)\wedge G(x))$ as $\forall_P x(G(x))$ and $\exists_P x(G(x))$, respectively; where $G$ is a formula. Finally, $\exists!x(G(x))$ is a metalinguistic abbreviation for $\exists x\forall y(G(y)\Leftrightarrow y = x)$.

Postulates {\sc F1}$\sim${\sc F11} of $\mbox{\boldmath{$\mathfrak{F}$}}$ are as follows.

\subsection{The very basic}

\begin{description}

\item[\sc F1 - Weak Extensionality] $\forall f \forall g (((f(g) = f \wedge g(f) = g) \vee (f(g) = g \wedge g(f) = f)) \Rightarrow f = g))$.

\end{description}

If $f(g) = f$ we say $f$ is {\em rigid with\/} $g$. If $f(g) = g$ we say $f$ is {\em flexible with\/} $g$. If both $f$ and $g$ are rigid (flexible) with each other, then $f = g$.

\begin{description}

\item[\sc F2 - Self-Reference] $\forall f (f(f) = f)$.

\end{description}

Any $f$ is flexible {\em and\/} rigid with itself. That may sound a strong limitation if we compare our framework, e.g., to lambda calculus \cite{Barendregt-13}. But that feature is quite useful for our purposes. One task for the future is to admit some variations of Flow where {\sc F2} is no axiom.

\begin{observacao}\label{Russell}

Let $y$ be a function such that $\forall x(y(x) = r \Leftrightarrow x(x)\neq r)$. What about $y(y)$? If $x = y$, then $y(y) = r \Leftrightarrow y(y)\neq r$ (Russell's paradox). But that entails $y(y)\neq y(y)$, another contradiction. However, our Self-Reference postulate does not allow us to define such an $y$, since $x(x) = x$ for any $x$. That is not the only way to avoid such an inconsistency. Our Restriction axiom (some pages below) does the same work. That opens the possibility of variations of Flow where Self-Reference is not a theorem.
\end{observacao}

Our first theorem states any $f$ can be identified by its images $f(x)$.

\begin{teorema}\label{igualdadefuncoes}
$\forall f\forall g(f = g\Leftrightarrow\forall x (f(x) = g(x)))$.
\end{teorema}

\begin{description}
\item[\sc Proof] From substitutivity of identity in $f(x) = f(x)$, proof of the $\Rightarrow$ part is straightforward: if $f = g$, then $f(x) = g(x)$, for any $x$. For the $\Leftarrow$ part, suppose, for any $x$, $f(x) = g(x)$. If $x = f$, $f(f) = g(f)$; if $x = g$, $f(g) = g(g)$. From {\sc F2}, $f(f) = f$ and $g(g) = g$. So, $g(f) = f$ and $f(g) = g$. {\sc F1} entails $f = g$.
\end{description}

Axioms {\sc F1} and {\sc F2} could be rewritten as one single formula:

\begin{description}

\item[\sc F1' - Alternative Weak Extensionality] $\forall f \forall g (((f(g) = f \wedge g(f) = g) \vee (f(g) = g \wedge g(f) = f)) \Leftrightarrow f = g))$.

\end{description}

{\sc F2} is a consequence from {\sc F1'}. Ultimately, $f=g$ entails $f(g) = f$ (from {\sc F1'}). And substitutivity of identity entails $f(f) = f$. But we prefer to keep axioms {\sc F1} and {\sc F2} (instead of {\sc F1'}) to smooth away some discussions. From {\sc F1} and {\sc F2}, we can analogously see that {\sc F1'} is a nontrivial theorem.

\begin{description}

\item[\sc F3 - Identity] $\exists f \forall x (f(x) = x)$.

\end{description}

There is at least one $f$ such that, for any $x$, we have $f(x) = x$. Any function $f$ which satisfies {\sc F3} is said to be an {\em identity function\/}.

\begin{teorema}\label{umbarraunico}
Identity function is unique.
\end{teorema}

\begin{description}
\item[\sc Proof] Suppose both $f$ and $g$ satisfy {\sc F3}. Then, for any $x$, $f(x) = x$ and $g(x) = x$. Thus, $f(g) = g$ and $g(f) = f$. Hence, according to {\sc F1}, $f = g$.
\end{description}

In other words, there is one single $f$ which is flexible to every term. In that case we simply say $f$ is {\em flexible\/}. That means ``flexible'' and ``identity'' are synonyms.

\begin{description}

\item[\sc F4 - Rigidness] $\exists f \forall x (f(x) = f)$.

\end{description}

There is at least one $f$ which is rigid with any function. Observe the symmetry between {\sc F3} and {\sc F4}! Any $f$ which satisfies last postulate is said to be {\em rigid\/}.

\begin{teorema}\label{unicidadezerobarra}
The rigid function is unique.
\end{teorema}

\begin{description}
\item[\sc Proof] Let $f$ and $g$ satisfy axiom {\sc F4}. Then, for any $x$, $f(x) = f$ and $g(x) = g$. Thus, $f(g(x)) = f(g) = f$ and $g(f(x)) = g(f) = g$. From {\sc F1}, $f = g$.
\end{description}

Now we are able to introduce new terminology. Term $\underline{\mathfrak{1}}$ is the identity (flexible) function, while $\underline{\mathfrak{0}}$ is the rigid function, since we proved they are both unique. So,
$$\forall x (\underline{\mathfrak{1}}(x) = x\;\;\wedge\;\;\underline{\mathfrak{0}}(x) = \underline{\mathfrak{0}}).$$
Since $f(x) = y$ says $f_1^2(f,x) = y$, {\sc F3} states there is an $f$ such that, for any $x$, $f_1^2(f,x) = x$, while {\sc F4} says there is $f$ where $f_1^2(f,x) = f$. If we do not grant the existence of other functions, it seems rather difficult to prove $\underline{\mathfrak{0}}\neq\underline{\mathfrak{1}}$.

\begin{teorema}\label{teoremamarcio1}
$\underline{\mathfrak{0}}$ is the only function which is rigid with $\underline{\mathfrak{0}}$.
\end{teorema}

\begin{description}
\item[\sc Proof] This theorem says $\forall x (x\neq \underline{\mathfrak{0}} \Rightarrow x(\underline{\mathfrak{0}}) \neq x)$, i.e., $\forall x (x(\underline{\mathfrak{0}}) = x \Rightarrow x = \underline{\mathfrak{0}})$. But $\underline{\mathfrak{0}}(x) = \underline{\mathfrak{0}}$. So, if  $x(\underline{\mathfrak{0}}) = x \wedge \underline{\mathfrak{0}}(x) = \underline{\mathfrak{0}}$, then $x = \underline{\mathfrak{0}}$ ({\sc F1}).
\end{description}

\begin{teorema}\label{teoremamarcio2}
$\underline{\mathfrak{1}}$ is the only function which is flexible with $\underline{\mathfrak{1}}$.
\end{teorema}

\begin{description}
\item[\sc Proof] This theorem says $\forall x (x\neq \underline{\mathfrak{1}} \Rightarrow x(\underline{\mathfrak{1}}) \neq \underline{\mathfrak{1}})$, i.e., $\forall x (x(\underline{\mathfrak{1}}) = \underline{\mathfrak{1}} \Rightarrow x = \underline{\mathfrak{1}})$. But $\underline{\mathfrak{1}}(x) = x$. So, if  $x(\underline{\mathfrak{1}}) = \underline{\mathfrak{1}} \wedge \underline{\mathfrak{1}}(x) = x$, then $x = \underline{\mathfrak{1}}$ ({\sc F1}).
\end{description}

For more details about $x(\underline{\mathfrak{0}})$ and $x(\underline{\mathfrak{1}})$, for any $x$, see Theorems \ref{fechamentorenato} and \ref{composicaodexcomumbarra}.

\begin{definicao}\label{actionaction}

$f[t]$ iff $t\neq f \wedge f(t) \neq \underline{\mathfrak{0}}$. We read $f[t]$ as `{\em $f$ acts on $t$\/}'.

\end{definicao}

Whereas $f(t)$ is a term, $f[t]$ abbreviates a formula. No $f$ acts on itself. The intuitive idea is to allow us to talk about what a function $f$ {\em does\/}. For example, $\underline{\mathfrak{0}}$ does nothing whatsoever, since there is no $t$ on which it acts.

Next definition is a pedagogical move to write postulate {\sc F5} in a simple and direct way. If the reader is not comfortable with that, all you have to do is to rewrite $h = f\circ g$ in {\sc F5} (next postulate) as the conjunction of all five formulas in the next definition.

\begin{definicao}\label{composicao}
Given $f$ and $g$, the {\em $\mathfrak{F}$-composition\/} $h = f\circ g$, if it exists, must satisfy the next conditions:

{\sc (i)} $h\neq\underline{\mathfrak{0}}$; {\sc (ii)} $\forall x((x\neq f \wedge x\neq g \wedge x\neq h)\Rightarrow h(x) = f(g(x)))$; \\{\sc (iii)} $h\neq f \Rightarrow h(f) = \underline{\mathfrak{0}}$; {\sc (iv)}  $h\neq g \Rightarrow h(g) = \underline{\mathfrak{0}}$; \\{\sc (v)} $(g\neq h\wedge f\neq h \wedge g\neq\underline{\mathfrak{1}}\wedge f\neq\underline{\mathfrak{1}})\Rightarrow (f(g(h)) = \underline{\mathfrak{0}}\vee g(h) = \underline{\mathfrak{0}})$

\end{definicao}

So, {\sc (i)} No $\mathfrak{F}$-composition is $\underline{\mathfrak{0}}$. {\sc (ii)} $\mathfrak{F}$-composition behaves like standard notions of composition between functions up to self-reference. Nevertheless, contrary to usual practice in standard set theories, we can define $\mathfrak{F}$-composition between any two functions. {\sc (iii and iv)} No $\mathfrak{F}$-composition $h = f\circ g$ acts on any of its factors $f$ or $g$. Finally, item {\sc (v)} is a technical constraint to avoid ambiguities in the calculation of $f\circ g$. Either $g$ does not act on the $\mathfrak{F}$-composition $h$ or, if it acts, then $f$ does not act on $g(h)$.

\begin{description}

\item[\sc F5 - $\mathfrak{F}$-Composition] $\forall f \forall g(\exists !h (h = f\circ g))$.

\end{description}

We never calculate $(f\circ g)(f\circ g)$ as $f(g(f\circ g))$, since $(f\circ g)(f\circ g)$ is $f\circ g$, according to {\sc F2}. Given $f$ and $g$, we can `build' $f\circ g$ through a three-step process: {\sc (i)} First we establish a label $h$ for $f\circ g$; {\sc (ii)} Next we evaluate $f(g(x))$ for any $x$ which is different of $f$ and $g$. By doing that we are assuming those $x$ are different of $h$. Thus, if such a choice of $x$ entails $f(g(x)) = y$ for a given $y$, then $h(x) = y$; {\sc (iii)} Next we evaluate the following possibilities: is $h$ equal to either $f$, $g$ or something else? If $h\neq g$, then $h(g)$ is supposed to be $\underline{\mathfrak{0}}$. If $h(g) = \underline{\mathfrak{0}}$ entails a contradiction, then $h$ is simply $g$. An analogous method is used for assessing if $h$ is $f$. Eventually, $h$ is neither $f$ nor $g$, as we can see in the next theorems.

It is worth to remark that $\mathfrak{F}$-composition plays a crucial role in our main result (Theorems \ref{PPvale} and \ref{AEfalha}).

\begin{teorema}\label{primeirophizero}

There is a unique $h$ such that $h\neq\underline{\mathfrak{0}}$ but $h(x) = \underline{\mathfrak{0}}$ for any $x\neq h$.

\end{teorema}

\begin{description}
\item[\sc Proof] From Definition \ref{composicao} and {\sc F5}, $\underline{\mathfrak{0}}\circ\underline{\mathfrak{0}}$ is a unique $h\neq\underline{\mathfrak{0}}$ such that, for any $x$ where $x\neq \underline{\mathfrak{0}}$ and $x\neq h$, $h(x) = \underline{\mathfrak{0}}(\underline{\mathfrak{0}}(x))=  \underline{\mathfrak{0}}(\underline{\mathfrak{0}}) = \underline{\mathfrak{0}}$. Since $h\neq\underline{\mathfrak{0}}$, then {\sc F5} entails $h(\underline{\mathfrak{0}}) = \underline{\mathfrak{0}}$. Thus, $h(x)$ is $\underline{\mathfrak{0}}$ for any $x\neq h$, while $h$ itself is different of $\underline{\mathfrak{0}}$.
\end{description}

Such $h$ of last theorem (which does not conflict neither with Theorem \ref{unicidadezerobarra} nor with Theorem \ref{teoremamarcio1}) is labeled with a special symbol, namely, $\varphi_0$. So, $\underline{\mathfrak{0}}\circ\underline{\mathfrak{0}} = \varphi_0$, where $\varphi_0\neq\underline{\mathfrak{0}}$. If the reader is intrigued by the subscript $_0$ in $\varphi_0$, our answer is `yes, we intend to introduce ordinals, where $\varphi_0$ is the first one'.

Let $f$, $g$, $a$, and $c$ be pairwise distinct functions where: {\sc (i)}$f(a) = g$, $f(g) = c$, $f(c) = c$, $f(f) = f$, and $f(r) = \underline{\mathfrak{0}}$ for the remaining values $r$; {\sc (ii)} $g(a) = a$, $g(g) = g$, and $g(r) = \underline{\mathfrak{0}}$ for the remaining values $r$; {\sc (iii)} $a$ and $c$ are arbitrary, as long they are neither $\underline{\mathfrak{0}}$ nor $\underline{\mathfrak{1}}$. Then, $(f\circ g)\circ f = \varphi_0$, while $f\circ (g\circ f) = l$, where $l(a) = c$, $l(l) = l$, and $l(r) = \underline{\mathfrak{0}}$ for the remaining values $r$. Thus, $\mathfrak{F}$-composition among functions which act on each other (observe in this example $f[g]$, since $f\neq g$ and $f(g) = c \neq \underline{\mathfrak{0}}$), if they exist, is not necessarily associative.

\begin{teorema}\label{associatividadeentrenaocirculares}
Let $f$, $g$, and $h$ be terms such that neither one of them acts on the remaining ones. Then $\mathfrak{F}$-Composition is associative among $f$, $g$, and $h$.
\end{teorema}

\begin{description}
\item[\sc Proof] Both $f\circ (g \circ h) = p$ and $(f\circ g) \circ h = q$ correspond (Definition \ref{composicao}) to $y = f(g(h(x)))$, if $x$ is different of $f$, $g$, $h$, $p$, and $q$. So, $f\circ (g \circ h) = (f\circ g) \circ h$ for those values of $x$. But $f$, $g$, and $h$ do not act on $f$, $g$, or $h$; and no $\mathfrak{F}$-composition acts on any of its factors. So, $f\circ (g \circ h) = (f\circ g) \circ h = f(g(h(x)))$.
\end{description}

Evaluation $f_1^2$ is not associative as well. Consider, e.g., $x(\underline{\mathfrak{1}}(x))$, for $x$ different of $\underline{\mathfrak{0}}$ and different of $\underline{\mathfrak{1}}$, and such that $x(\underline{\mathfrak{1}}) = \underline{\mathfrak{0}}$ (functions like that do exist, as we can see later on). Thus, $x(\underline{\mathfrak{1}}(x)) = x(x) = x$. If evaluation was associative, we would have $x(\underline{\mathfrak{1}}(x)) = x(\underline{\mathfrak{1}})(x) = \underline{\mathfrak{0}}(x) = \underline{\mathfrak{0}}$; a contradiction, since we assumed $x\neq\underline{\mathfrak{0}}$! Of course this {\em rationale\/} works only if we prove the existence of a function like $x$ and that $\underline{\mathfrak{0}}\neq\underline{\mathfrak{1}}$. Both claims are the subject of the next two theorems.

It is good news that evaluation is not associative. According to {\sc F2}, for all $t$, $g(t) = (g(g))(t)$, since $g(g) = g$. If evaluation was associative, we would have $g(t) = g(g(t))$ and, thus, $g = g\circ g$. So, $\mathfrak{F}$-composition would be idempotent, a useless feature for our purposes.

\begin{teorema}

$\underline{\mathfrak{0}}\neq\underline{\mathfrak{1}}$.

\end{teorema}

\begin{description}
\item[\sc Proof] $\underline{\mathfrak{0}}\circ\underline{\mathfrak{0}} = \varphi_0\neq\underline{\mathfrak{0}}$ (Theorem \ref{primeirophizero}). But $\underline{\mathfrak{0}}(\varphi_0) = \underline{\mathfrak{0}}$, while $\underline{\mathfrak{1}}(\varphi_0) = \varphi_0$. Thus, from Theorem \ref{igualdadefuncoes}, $\underline{\mathfrak{0}}\neq\underline{\mathfrak{1}}$.
\end{description}

\begin{teorema}\label{segundophizero}

For any $x$ we have $\underline{\mathfrak{0}}\circ x = x\circ\underline{\mathfrak{0}} = \varphi_0$

\end{teorema}

\begin{description}
\item[\sc Proof] First we prove $\underline{\mathfrak{0}}\circ x = \varphi_0$. From Definition \ref{composicao} and axiom {\sc F5}, $(\underline{\mathfrak{0}}\circ x)(t) = \underline{\mathfrak{0}}(x(t)) = \underline{\mathfrak{0}}$ for any $t\neq \underline{\mathfrak{0}}$ and $t\neq x$. But {\sc F5} demands $\underline{\mathfrak{0}}\circ x$ is different of $\underline{\mathfrak{0}}$. Hence, $(\underline{\mathfrak{0}}\circ x)(\underline{\mathfrak{0}}) = \underline{\mathfrak{0}}$. Regarding $x$, there are three possibilities: ({\sc i}) $x = \underline{\mathfrak{0}}$; ({\sc ii}) $x = \varphi_0$; ({\sc iii}) $x$ is neither $\underline{\mathfrak{0}}$ nor $\varphi_0$. The first case corresponds to Theorem \ref{primeirophizero}, which entails $\underline{\mathfrak{0}}\circ x = \varphi_0$. In the second case, if $\underline{\mathfrak{0}}\circ x$ is different of $x = \varphi_0$, then $(\underline{\mathfrak{0}}\circ x)(\varphi_0) = \underline{\mathfrak{0}}$. But that would entail $(\underline{\mathfrak{0}}\circ x)(t) = \underline{\mathfrak{0}}$ for any $t\neq \underline{\mathfrak{0}}\circ x$, which corresponds exactly to function $\varphi_0$ proven in Theorem \ref{primeirophizero}, a contradiction. So, $\underline{\mathfrak{0}}\circ x$ is indeed $\varphi_0$, when $x = \varphi_0$. Concerning last case, since $x\neq\underline{\mathfrak{0}}$ and $x\neq\varphi_0$, then (Theorem \ref{igualdadefuncoes}) there is $t\neq\varphi_0$ such that $x(t)\neq\underline{\mathfrak{0}}$ for $x\neq t$. From {\sc F5}, $(\underline{\mathfrak{0}}\circ x)(t) = \underline{\mathfrak{0}}$ for such value of $t$. But once again we have a function $\underline{\mathfrak{0}}\circ x$ such that $(\underline{\mathfrak{0}}\circ x)(t) = \underline{\mathfrak{0}}$ for any $t\neq \underline{\mathfrak{0}}$, which corresponds to $\varphi_0$ from Theorem \ref{primeirophizero}. Concerning the identity $x\circ\underline{\mathfrak{0}} = \varphi_0$, the proof is analogous. If $h = x\circ\underline{\mathfrak{0}}$, then $h\neq\underline{\mathfrak{0}}$ (from {\sc F5}). Besides, for any $x$, $x\neq h$ entails $h(x) = \underline{\mathfrak{0}}$. Hence, $h = \varphi_0$.
\end{description}

\begin{teorema}\label{fechamentorenato}

$\forall x(x(\underline{\mathfrak{0}}) = \underline{\mathfrak{0}})$.

\end{teorema}

\begin{description}
\item[\sc Proof] $x\circ \underline{\mathfrak{0}} = \varphi_0$ (Theorem \ref{segundophizero}). Hence, for any $t$, $(t\neq x\wedge t\neq\underline{\mathfrak{0}}\wedge t\neq\varphi_0)\Rightarrow \varphi_0(t) = x(\underline{\mathfrak{0}}(t)) = x(\underline{\mathfrak{0}})$. But $\varphi_0(t) = \underline{\mathfrak{0}}$ for any $t\neq\varphi_0$. Thus, $x(\underline{\mathfrak{0}}) = \underline{\mathfrak{0}}$.
\end{description}

\begin{teorema}\label{umbarracompostocomumbarraehumbarra}

$\underline{\mathfrak{1}}\circ\underline{\mathfrak{1}} = \underline{\mathfrak{1}}$.

\end{teorema}

\begin{description}
\item[\sc Proof]
From {\sc F5}, $h = \underline{\mathfrak{1}}\circ\underline{\mathfrak{1}}$ entails $h(x) = x$ for any $x\neq\underline{\mathfrak{1}}$. Since the $\mathfrak{F}$-composition is unique and $\underline{\mathfrak{1}}$ guarantees all demanded conditions, then $h = \underline{\mathfrak{1}}$.
\end{description}

\begin{teorema}\label{dezenovedemarco}

$\forall x \forall y(x\circ y = \underline{\mathfrak{1}}\Rightarrow (x = \underline{\mathfrak{1}} \wedge y = \underline{\mathfrak{1}}))$.

\end{teorema}

\begin{description}
\item[\sc Proof] Let $x\neq\underline{\mathfrak{1}}$. Then $\underline{\mathfrak{1}}(x) = \underline{\mathfrak{0}}$, according to {\sc F5}. But that happens only for $x = \underline{\mathfrak{0}}$. And $\underline{\mathfrak{0}}\circ y\neq\underline{\mathfrak{1}}$. Analogous argument holds for $y\neq\underline{\mathfrak{1}}$.
\end{description}

There can be no functions different of $\underline{\mathfrak{1}}$ such that their composition is $\underline{\mathfrak{1}}$.

Next theorem is important for a better understanding about {\sc F1}, although its proof does not demand the use of such a postulate.

\begin{teorema}\label{composicaodexcomumbarra}

$\forall x((x\neq \underline{\mathfrak{0}} \wedge x(\underline{\mathfrak{1}}) = \underline{\mathfrak{0}})\Leftrightarrow (\underline{\mathfrak{1}}\circ x = x\wedge x\circ\underline{\mathfrak{1}} = x \wedge x\neq\underline{\mathfrak{1}}))$.

\end{teorema}

\begin{description}
\item[\sc Proof] The $\Rightarrow$ part. If $x(\underline{\mathfrak{1}}) = \underline{\mathfrak{0}}$, then $x\neq\underline{\mathfrak{1}}$, since $\underline{\mathfrak{1}}(\underline{\mathfrak{1}}) = \underline{\mathfrak{1}}$. If $x\circ \underline{\mathfrak{1}} = h$, then, for any $t$ different of $x$, $\underline{\mathfrak{1}}$, and $h$, we have $h(t) = x(\underline{\mathfrak{1}}(t)) = x(t)$. If $h = x$, then $h$ satisfies all conditions from {\sc F5}, since $x\neq\underline{\mathfrak{1}}$, $h(\underline{\mathfrak{1}}) = x(\underline{\mathfrak{1}}) = \underline{\mathfrak{0}}$, and $h(t) = x(t)$ for any $t\neq h$. Since {\sc F5} demands $h$ to be unique, then $h = x$. If $\underline{\mathfrak{1}}\circ x = h$, we use an analogous argument. For the $\Leftarrow$ part, $\underline{\mathfrak{1}}\circ x = x\circ\underline{\mathfrak{1}} = x$ entails $x\neq\underline{\mathfrak{0}}$ (Theorem \ref{segundophizero}). Since $x\neq\underline{\mathfrak{1}}$, then {\sc F5} demands for the $\mathfrak{F}$-composition $x$ that $x(\underline{\mathfrak{1}}) = \underline{\mathfrak{0}}$.
\end{description}

\begin{observacao}\label{observaunicidade}

If it wasn't for the uniqueness of $\mathfrak{F}$-compositions in {\sc F5}, Flow would be consistent with the existence of many functions, like $x$ and $h$ (from the proof of Theorem \ref{composicaodexcomumbarra}), which ``do'' the same thing. We refer to such functions as {\em clones\/}. For a brief investigation about clones see Definition \ref{definesim} and its subsequent discussion. Clones are meant to be different functions $x$ and $h$ which share the same images $x(t)$ and $h(t)$ for any $t$ different of both $x$ and $h$, and such that $x(h) = \underline{\mathfrak{0}}$ and $h(x) = \underline{\mathfrak{0}}$. From Theorem \ref{igualdadefuncoes}, $x\neq h$ (recall {\sc F2}). We use this opportunity to prove Theorem \ref{primeirophizero}, since $\underline{\mathfrak{0}}\circ\underline{\mathfrak{0}} = \varphi_0$, where $\varphi_0$ and $\underline{\mathfrak{0}}$ are clones. But we cease using clones when we talk about other functions. That is why we refer to {\sc F1} as ``weak extensionality''. A strong extensionality postulate would demand that other clones besides $\underline{\mathfrak{0}}$ and $\varphi_0$ cannot exist. We do something like that in some remaining postulates where we use the quantifier $\exists!$.
\end{observacao}

\begin{definicao}\label{defineSigmasigma}

Given a term $f$, $g$ is the $\mathfrak{F}$-successor of $f$, and we denote this by $\Sigma(f,g)$, iff $f(g) = \underline{\mathfrak{0}} \wedge \forall x(x\neq g \Rightarrow g(x) = f(x))$.

\end{definicao}

For example, let $f$ be given by $f(a) = b$, $f(b) = c$, $f(c) = c$, $f(f) = f$, and $f(r) = \underline{\mathfrak{0}}$ for the remaining values $r$. Now let $g$ be given by $g(a) = b$, $g(b) = c$, $g(c) = c$, $g(f) = f$, $g(g) = g$, and $g(r) = \underline{\mathfrak{0}}$ for the remaining values $r$, where $g\neq f$. In that case, $\Sigma(f,g)$, if functions like those exist. Observe this has nothing to do with our previous discussion about clones.

\begin{description}

\item[\sc F6 - $\mathfrak{F}$-Successor Function] $\exists!\sigma(\sigma\neq\underline{\mathfrak{0}}\wedge \forall f((f\neq\sigma\wedge f\neq\underline{\mathfrak{0}})\Rightarrow \\(\exists g(\Sigma(f,g)\Leftrightarrow \sigma(f) = g)\vee (\forall h(\neg\Sigma(f,h))\Leftrightarrow \sigma(f) = \underline{\mathfrak{0}}))))$.

\end{description}

This last axiom states the existence and uniqueness of a special function $\sigma$. From now on we write mostly $\sigma_f$ for $\sigma(f)$. Every time we use the symbol $\sigma$ we are referring to the same term from {\sc F6}: the only one such that $\sigma_f = g \neq \underline{\mathfrak{0}}$ is equivalent to $\Sigma(f,g)$, and $\sigma_f = \underline{\mathfrak{0}}$ is equivalent to $\forall h(\neg\Sigma(f,h))$, as long $f$ is neither $\sigma$ nor $\underline{\mathfrak{0}}$. Thus, $\sigma$ successfully `signals' $\mathfrak{F}$-successors, when they exist, for all functions, except when $f$ is either $\underline{\mathfrak{0}}$ or $\sigma$. In those cases, we have $\Sigma(\underline{\mathfrak{0}},\varphi_0)$ and $\sigma_{\underline{\mathfrak{0}}} = \underline{\mathfrak{0}}$ (Theorem \ref{fechamentorenato}), while $\forall h (\neg\Sigma(\sigma, h))$ and $\sigma_{\sigma} = \sigma$ ({\sc F2}).

Next we want to grant the existence of $\mathfrak{F}$-successors for many other terms. As previously announced, this axiom states the existence of an hierarchy of functions, where $\varphi_0$ is the first one.

\begin{description}

\item[\sc F7 - Infinity] $\exists i ((\forall t (i(t) = t \vee i(t) = \underline{\mathfrak{0}}))\wedge \sigma_i \neq \underline{\mathfrak{0}} \wedge (i(\varphi_0) = \varphi_0 \wedge \forall x ((x\neq\underline{\mathfrak{0}}\wedge i(x) = x) \Rightarrow (i(\sigma_x) = \sigma_x \wedge \sigma_x\neq\underline{\mathfrak{0}}))))$.

\end{description}

\begin{definicao}\label{funcaoindutiva}
Any $i$ which satisfies {\sc F7} is said to be {\em inductive\/}.
\end{definicao}

Since the existence of $\varphi_0$ is granted by {\sc F5}, we can use $\sigma$ and {\sc F7} to grant the existence of $\varphi_1 = \sigma_{\varphi_0}$ such that $\varphi_1(\varphi_1) = \varphi_1$, $\varphi_1(\varphi_0) = \varphi_0$, and for the remaining values $r$ (those who are neither $\varphi_0$ nor $\varphi_1$) we have $\varphi_1(r) = \underline{\mathfrak{0}}$. That happens because {\sc F7} states the existence of at least one other function $i$ and infinitely many other functions. It says $i(\varphi_0) = \varphi_0$. Besides, there is a non-$\underline{\mathfrak{0}}$ $\mathfrak{F}$-successor of $\varphi_0$ such that $i(\sigma_{\varphi_0}) = \sigma_{\varphi_0}$. More than that, if $x$ admits a non-$\underline{\mathfrak{0}}$ $\mathfrak{F}$-successor $\sigma_x$ (where $i(x) = x$), then $i(\sigma_x) = \sigma_x$, where $\sigma_x \neq \underline{\mathfrak{0}}$. Thus, $\varphi_0\neq\underline{\mathfrak{0}}$ and $\sigma_{\varphi_0} = \varphi_1$, where $\varphi_1\neq\varphi_0$ and $\varphi_1\neq\underline{\mathfrak{0}}$. Analogously we can get $\varphi_2$, $\varphi_3$, and so on. Along with those terms $\varphi_n$, {\sc F7} says any inductive function $i$ admits its own non-$\underline{\mathfrak{0}}$ $\mathfrak{F}$-successor $\sigma_i$.

Subscripts $0$, $1$, $2$, $3$, etc., are metalinguistic symbols based on an alphabet of ten symbols (the usual decimal numeral system) which follows the lexicographic order $\prec$, where $0\prec 1\prec 2\prec \cdots\prec 8\prec 9$. If $n$ is a subscript, then $n+1$ corresponds to the next subscript, in accordance to the lexicographic order. In that case, we write $n\prec n+1$. $n+m$ is an abbreviation for $(...(...((n+1)+1)+...1)...)$ with $m$ occurrences of $+$ and $m$ occurrences of pairs of parentheses. Again we have $n\prec n+m$. Besides, $\prec$ is a strict total order. That fact allows us to talk about a minimum value between subscripts $m$ and $n$: $\mbox{min}\{ m, n\}$ is $m$ iff $m\prec n$, it is $n$ iff $n\prec m$, and it is either one of them if $m = n$. Of course, $m = n$ iff $\neg(m\prec n) \wedge \neg (n\prec m)$. If $m\prec n \vee m = n$, we denote this by $m\preceq n$. Such a vocabulary of ten symbols endowed with $\prec$ is called here {\em (meta) language\/} $\mathcal L$.

{\sc F7} provides us some sort of ``recursive definition'' for functions $\varphi_n$, while it allows as well to guarantee the existence of inductive functions: {\sc (i)} $\varphi_0$ is such that $\varphi_0(x)$ is $\varphi_0$ if $x = \varphi_0$ and $\underline{\mathfrak{0}}$ otherwise; {\sc (ii)} $\varphi_{n+1}$ is such that $\varphi_{n+1}(\varphi_{n+1}) = \varphi_{n+1}$, $\varphi_{n+1}\neq \varphi_n$, and $\varphi_{n+1}(x) = \varphi_n(x)$ for any $x$ different of $\varphi_{n+1}$.

Observe that $\varphi_{n+1}(\varphi_n) = \varphi_n(\varphi_n) = \varphi_n$, while $\varphi_n(\varphi_{n+1}) = \underline{\mathfrak{0}}$. Moreover, $\varphi_{n+2}(\varphi_{n+1}) = \varphi_{n+1}$, and $\varphi_{n+2}(\varphi_n) = \varphi_{n+1}(\varphi_n) = \varphi_n$; while $\varphi_n(\varphi_{n+2}) = \underline{\mathfrak{0}}$.

{\sc Figure 1} illustrates how to represent some functions $f$ in a quite intuitive way. A {\em diagram\/} of $f$ is formed by a rectangle. On the left top corner inside the rectangle we find label $f$. The remaining labels refer to terms $x$ such that either $f(x)\neq \underline{\mathfrak{0}}$ ($f$ acts on $x$) or there is $t$ such that $f(t) = x \neq \underline{\mathfrak{0}}$. Term $\underline{\mathfrak{0}}$ never occurs in any diagram. For each $x$ inside the rectangle there is a unique corresponding arrow which indicates the image of $x$ by $f$, as long $x$ is not $f$ itself. Due to self-reference, label $f$ at the left top corner inside the rectangle does not need to be attached to any arrow, to avoid redundancy.

From left to right, the first diagram in {\sc Figure 1} refers to $\varphi_0$. It says, for any $x$, $\varphi_0(x)$ is $\underline{\mathfrak{0}}$, except for $\varphi_0$ itself. The second diagram says $\varphi_1(\varphi_1) = \varphi_1$, and $\varphi_1(\varphi_0) = \varphi_0$. The circular arrow associated to $\varphi_0$ in the second diagram says $\varphi_1(\varphi_0) = \varphi_0$ ($\varphi_1[\varphi_0]$). The third diagram says $\varphi_2(\varphi_2) = \varphi_2$, $\varphi_2(\varphi_1) = \varphi_1$, and $\varphi_2(\varphi_0) = \varphi_0$. Finally, for the sake of illustration, the last diagram corresponds to an $f$ such that $f(a) = b$, $f(b) = c$, and $f(c) = \underline{\mathfrak{0}}$. That is why there is no arrow `starting' at $c$. In those cases we represent $c$ outside the rectangle. The existence of functions like that is granted by {\sc F10}$_{\alpha}$ some pages below. The diagrams of $\underline{\mathfrak{0}}$ and $\underline{\mathfrak{1}}$ are, respectively, a blank rectangle and a filled in black rectangle. Other examples are provided in the next paragraphs.\\

\begin{picture}(60,50)

\put(3,13){\framebox(50,40)}


\put(10,43){$\varphi_0$}


\put(63,13){\framebox(50,40)}


\put(70,43){$\varphi_1$}

\put(70,26){$\boldmath \curvearrowleft$}

\put(70,20){$\varphi_0$}


\put(123,13){\framebox(50,40)}


\put(130,43){$\varphi_2$}

\put(130,26){$\boldmath \curvearrowleft$}

\put(130,20){$\varphi_0$}

\put(160,26){$\boldmath \curvearrowleft$}

\put(160,20){$\varphi_1$}

\put(190,30){$\cdots$}

\put(250,13){\framebox(50,40)}

\put(255,43){$f$}

\put(255,23){$a$}

\put(275,23){$b$}

\put(315,23){$c$}

\put(282,26){\vector(2,0){30}}

\put(262,26){\vector(2,0){10}}

\end{picture}

\vspace{-3mm}

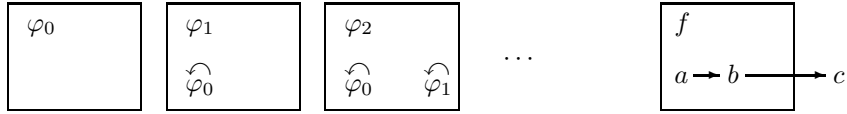
{\sc Figure 1:} From left to right, diagrams of $\varphi_0$, $\varphi_1$, $\varphi_2$, and an arbitrary $f$.

Observe that $\varphi_{m+n}(\varphi_n) = \varphi_n$, for any $m$ and $n$ of $\mathcal L$.

\begin{teorema}\label{noaction}

$\underline{\mathfrak{0}}$ and $\varphi_0$ are the only functions who do not act on any $t$.

\end{teorema}

\begin{description}
\item[\sc Proof] If $f$ does not act on any $t$, then $\forall t (t = f\vee f(t) = \underline{\mathfrak{0}})$. Suppose $f$ is neither $\underline{\mathfrak{0}}$ nor $\varphi_0$. Then there is $t$ such that $t\neq f \wedge f(t)\neq\underline{\mathfrak{0}}$. But that entails $f[t]$ (Definition \ref{actionaction}). So, the only functions which do not act on any $t$ are $\underline{\mathfrak{0}}$ and $\varphi_0$.
\end{description}

\begin{teorema}\label{sucessordeumbarra}
$\sigma_{\underline{\mathfrak{1}}} = \underline{\mathfrak{0}}$.
\end{teorema}

\begin{description}
\item[\sc Proof] $\sigma_{\underline{\mathfrak{1}}}\neq\underline{\mathfrak{0}}\Rightarrow\underline{\mathfrak{1}}(\sigma_{\underline{\mathfrak{1}}}) = \underline{\mathfrak{0}}$ (Definition \ref{defineSigmasigma}). That happens only if $\sigma_{\underline{\mathfrak{1}}} = \underline{\mathfrak{0}}$.
\end{description}

\begin{teorema}\label{psi}

For any $g$, $\sigma_g \neq \underline{\mathfrak{1}}$

\end{teorema}

\begin{description}
\item[\sc Proof] Suppose there is $g$ such that $\sigma_g = \underline{\mathfrak{1}}$. Since $\underline{\mathfrak{1}}\neq\underline{\mathfrak{0}}$, then $g(\underline{\mathfrak{1}}) = \underline{\mathfrak{0}}$ and $\underline{\mathfrak{1}}(g) = g$, according to Definition \ref{defineSigmasigma}. So, $g\neq\underline{\mathfrak{1}}$. Once again from Definition \ref{defineSigmasigma}, $g(x) = \underline{\mathfrak{1}}(x)$ for any $x\neq\underline{\mathfrak{1}}$. But those are the same conditions for $\underline{\mathfrak{1}}\circ\underline{\mathfrak{1}} = g$, according to {\sc F5}. Since any $\mathfrak{F}$-composition is unique, then $g = \underline{\mathfrak{1}}$ (Theorem \ref{umbarracompostocomumbarraehumbarra}), which contradicts the assumption $g\neq\underline{\mathfrak{1}}$. So, there is no such $g$.
\end{description}

In other words, the existence of some `functions' in Flow is forbidden.

\begin{definicao}\label{restricaodefuncao}
For any function $f$, a {\em restriction\/} $g$ {\em of\/} $f$ is defined as

\noindent
$$g\subseteq f \;\;\mbox{iff}\;\; g\neq\underline{\mathfrak{0}}\wedge\forall x((g[x]\Rightarrow f[x]) \wedge ((g[x]\wedge f[x])\Rightarrow f(x) = g(x))),$$

\end{definicao}

Proper restrictions are defined as $g\subset f\;\;\mbox{iff}\;\; g\subseteq f \wedge g\neq f$. We abbreviate $\neg (g\subseteq f)$ and $\neg (g\subset f)$ as, respectively, $g\not\subseteq f$ and $g\not\subset f$. For example, $\varphi_1\subseteq \varphi_3$, $\varphi_1\subset \varphi_3$, and $\varphi_3\not\subseteq \varphi_1$.

\begin{teorema}\label{restricaoderenato}

$\forall f \forall g (g\subset f \Rightarrow g(f) = \underline{\mathfrak{0}})$.

\end{teorema}

\begin{description}
\item[\sc Proof] Suppose $g(f)\neq\underline{\mathfrak{0}}$. Since $g\neq f$, then $g[f]$. But for any $f$ we have $\neg f[f]$. In other words, $\neg (g[f]\Rightarrow f[f])$. Hence, $g\not\subseteq f$.
\end{description}

Some of the most useful restrictions are obtained from a given formula $F$, in a way which resembles the well known Separation Scheme in ZF. That is achieved thanks to careful considerations regarding the $\mathfrak{F}$-successor function $\sigma$. But before that, we need more concepts, since we are interested on a vast number of situations.

\begin{definicao}\label{abrangente}

A term $f$ is {\em comprehensive\/} iff there is $g$ such that $g\neq\underline{\mathfrak{0}}$, $g\subseteq f$, and $\sigma_g = \underline{\mathfrak{0}}$; and we denote that by $\mathbb{C}(f)$. Otherwise, $f$ is {\em uncomprehensive\/}.

\end{definicao}

Comprehensive functions are supposed to describe ``huge'' functions who act on ``many terms'', and they are partially regulated by some of the the next axioms. If $f$ itself is such that $\sigma_f = \underline{\mathfrak{0}}$, then $f$ is comprehensive, except when $f = \underline{\mathfrak{0}}$. Besides, it is easy to see that $\underline{\mathfrak{1}}$ is comprehensive and $\underline{\mathfrak{0}}$ is uncomprehensive. On the other hand, there are other comprehensive functions which, by the way, play an important role within our proposal for building a model of ZFU. But for now we are mostly interested on uncomprehensive functions, as it follows in the next Subsection.

\begin{teorema}\label{extensaodeabrangenteehabrangente}
If $f\subseteq g$ and $f$ is comprehensive, then $g$ is comprehensive.
\end{teorema}

\begin{description}
\item[\sc Proof] If $f$ is comprehensive, then there is $h$ such that $h\neq\underline{\mathfrak{0}}$, $h\subseteq f$, and $\sigma_h=\underline{\mathfrak{0}}$. Since $h\subseteq f\subseteq g$, then $h\subseteq g$. Thus, $g$ is comprehensive.
\end{description}

\subsection{Emergent functions}

\begin{definicao}\label{definindoemergente}

$\mathbb{E}(f)$ iff {\sc (i)} $f\neq\sigma$; {\sc (ii)} $\sigma_f\neq\underline{\mathfrak{0}}$; {\sc (iii)} $\forall x(f[x]\Rightarrow \sigma_x\neq\underline{\mathfrak{0}})$; {\sc (iv)} $\forall y((\exists x(f(x) = y \wedge y\neq\underline{\mathfrak{0}}))\Rightarrow \sigma_y\neq\underline{\mathfrak{0}})$.

\end{definicao}

We read $\mathbb{E}(f)$ as `$f$ is {\em emergent\/}'. $\varphi_0$ is vacuously emergent. That entails $\varphi_1$ is emergent. Emergent functions are supposed to be uncomprehensive terms who act only on uncomprehensive terms. The behavior of emergent functions is regulated by the next axiom, {\sc F8}.

\begin{definicao}\label{lurking}

$g\trianglelefteq f$ iff $g\neq\underline{\mathfrak{0}}\wedge\forall x(g[x]\Rightarrow ((f[x]\vee \exists a(f(a) = x))\wedge (f[g(x)]\vee \exists b (f(b) = g(x)))))$. If $g\trianglelefteq f$, we say $g$ {\em lurks\/} $f$. Besides, $g\triangleleft f$ iff $g\trianglelefteq f$ and $g\neq f$. In that case we say $g$ {\em properly lurks\/} $f$. Finally, $\neg (g\trianglelefteq f)$ and $\neg (g\triangleleft f)$ are abbreviated as $g\ntrianglelefteq f$ and $g\ntriangleleft f$, respectively.

\end{definicao}

Last definition is a generalization of restriction, as it follows in the next theorem. Besides, if $f$ is a function such that $f(f) = f$, $f(\varphi_0) = \varphi_1$, $f(\varphi_1) = \varphi_0$ and $f(r) = \underline{\mathfrak{0}}$ for the remaining values, then $f\trianglelefteq \varphi_2$, $\varphi_2\trianglelefteq f$, $f\trianglelefteq \varphi_3$, $\varphi_3\ntrianglelefteq f$.

\begin{teorema}

If $g\subseteq f$ and $g\neq\underline{\mathfrak{0}}$, then $g\trianglelefteq f$.

\end{teorema}

The proof is straightforward. The converse is obviously not a theorem.

\begin{definicao}\label{definindomaximapotencia}

$h = \mathfrak{p}(f)$ iff $\forall x (h[x]\Leftrightarrow x\trianglelefteq f)$. $h$ is the {\em full power\/} of $f$.

\end{definicao}

The full power $\mathfrak{p}(f)$ of $f$ acts on all terms $x$ who lurk $f$.

\begin{teorema}

If $f\subset\underline{\mathfrak{1}}$ acts on $n$ terms, then $\mathfrak{p}(f)$ acts on $(n+1)^n$ terms.

\end{teorema}

\begin{description}
\item[\sc Proof] If $f$ acts on $n$ terms $x_i$ and $g$ lurks $f$, each $x_i$ may correspond to any $x_j$ (where eventually $x_j = x_i$) in the sense we may have $g(x_i) = x_j$. On the other hand, we may have $g(x_i) = \underline{\mathfrak{0}}$ as well. Thus, for each $x_i$ ($n$ possible values) there are $n+1$ possible images. So, there are $(n+1)^n$ terms who lurk $f$.
\end{description}

\begin{description}

\item[\sc F8 - Cohesion] $\forall f(\mathbb{E}(f)\Rightarrow (\forall g ((g\trianglelefteq f\Rightarrow \sigma_g\neq\underline{\mathfrak{0}})\wedge \forall h ((h[g]\Leftrightarrow g\trianglelefteq f)\Rightarrow \sigma_h\neq\underline{\mathfrak{0}}))\wedge \forall i(\forall t(i[t]\Rightarrow \exists x(x[t]\wedge f[x])))\Rightarrow \sigma_i\neq\underline{\mathfrak{0}})\wedge \mathbb{E}(\sigma_f))$.

\end{description}

If $f$ is emergent, then any $g$ (if it exists) who lurks $f$ has a non-$\underline{\mathfrak{0}}$ $\mathfrak{F}$-successor. Function $h$ at {\sc F8} refers to the full power of $f$, if it exists. Term $i$ is useful for dealing with arbitrary union, if we are able to define it. Finally, if $f$ is emergent, so it is $\sigma_f$.

Observe {\sc F8} is an existence postulate. For understanding this, recall {\sc F6}, which provides necessary and sufficient conditions for knowing if $\sigma_f\neq\underline{\mathfrak{0}}$: {\em there must be\/} a $g$ different of $f$ such that certain conditions are met. The point here is that {\sc F8} grants the existence of certain terms which are $\mathfrak{F}$-successors of others, as long some conditions are met. When we say, as above, that $\sigma_g\neq\underline{\mathfrak{0}}$, we state {\em there is\/} a function $z\neq\underline{\mathfrak{0}}$ such that $\sigma_g = z$. The same happens to $\sigma_h$ and $\sigma_i$.

\subsection{More about restrictions}

\begin{definicao}\label{delimitadorarestricao}

Let $f$ be a function and $F(t)$ be a formula where all occurrences of $t$ are free. We say $g$ is {\em restriction of $f$ under\/} $F(t)$, and denote this by $g = f\big|_{F(t)}$, iff: {\sc (i)} $g\neq\underline{\mathfrak{0}}$; {\sc (ii)} $f(\sigma_g) = \underline{\mathfrak{0}} \vee \neg F(g)$; {\sc (iii)} $g\neq f\Rightarrow g(f) = \underline{\mathfrak{0}}$; {\sc (iv)} $\forall t((t\neq f\wedge t\neq g)\Rightarrow ((g(t) = f(t)\wedge F(t)\wedge f[t])\vee (g(t) = \underline{\mathfrak{0}}\wedge (\neg F(t)\vee \neg f[t]))))$.

\end{definicao}

Formula $g = f\big|_{F(t)}$ is somehow equivalent to $g\subseteq f$, as we see in the next two theorems.

\begin{teorema}

If $g\subseteq f$, then it is possible to state a formula $F$ where $g = f\big|_F$.

\end{teorema}

\begin{description}
\item[\sc Proof] If $g\subseteq f$, assume as formula $F(t)$ the next one: $g[t]$. Item ({\sc i}) of Definition \ref{delimitadorarestricao} is a consequence from Definition \ref{restricaodefuncao}. Item ({\sc ii}) of the same definition is granted thanks to the fact that $\neg g[g]$ ($\neg F(g)$) for any $g$. Item ({\sc iii}) is due to Theorem \ref{restricaoderenato}. Finally, item ({\sc iv}) is granted from Definition \ref{restricaodefuncao}.
\end{description}

\begin{teorema}

If $g = f\big|_F$, then $g\subseteq f$.

\end{teorema}

\begin{description}
\item[\sc Proof] If $g = f\big|_F$, then either $g = \varphi_0$ or $g\neq\varphi_0$. In the first case, the proof is immediate by vacuity ($\varphi_0$ does not act on any term). If $g\neq\varphi_0$, then item ({\sc iv}) of Definition \ref{delimitadorarestricao} demands, for any $t$, $g[t]\Rightarrow (f[t]\wedge g(t) = f(t))$. Thus, $g\subseteq f$.
\end{description}

A natural way of getting some restrictions $g$ of $f$ is through {\sc F9}$_{\mathbb{E}}$ below. Observe as well item ({\sc iv}) of Definition \ref{delimitadorarestricao} takes into account the self-reference postulate, since we demand $t\neq f \wedge t\neq g$. Now, if $F(t)$ is a formula (abbreviated by $F$) where all occurrences of $t$ are free, then the following is an axiom.

\begin{description}

\item[\sc F9$_{\mathbb{E}}$ - $\mathbb{E}$-Restriction] $\forall f\forall x((F(x)\Rightarrow \mathbb{E}(x))\Rightarrow \exists! g(g = f\big|_F))$.

\end{description}

Subscript $_{\mathbb{E}}$ highlights a strong commitment to emergent functions, although there is no need of $f$ to be emergent. Many restrictions due to this last postulate grant the existence of emergent functions. Next theorem, for example, shows we do not need $\mathfrak{F}$-composition to prove there is $\varphi_0$.

\begin{teorema}\label{teoremasobrerestricaodezerobarra}

$\underline{\mathfrak{0}}\big|_{F(x)} = \varphi_0$ if $\forall x(F(x)\Rightarrow \mathbb{E}(x))$.

\end{teorema}

\begin{description}
\item[\sc Proof] Immediate, since $\underline{\mathfrak{0}}$ does not act on any term and no restriction can be $\underline{\mathfrak{0}}$.
\end{description}

\begin{teorema}\label{trioderestricoes}

For any emergent $f$ we have: {\sc (i)} $f\big|_{x\neq x} = \varphi_0$; {\sc (ii)} $f\big|_{x = x} = f$; {\sc (iii)} $f\big|_{x \neq f} = f$; {\sc (iv)} $f\big|_{x = f} = \varphi_0$.

\end{teorema}

\begin{description}
\item[\sc Proof]Item {\sc (i)} is proven by vacuity. {\sc (ii)} takes into account all terms where $f$ acts. Observe {\sc F8} grants any emergent function who acts on any term, acts on emergent functions. About {\sc (iii)} and {\sc (iv)}, recall $f$ plays no role into the calculation of its restriction. All that matters are the terms where $f$ acts.
\end{description}

From {\sc F9}$_{\mathbb{E}}$, there are {\em four\/} possible restrictions $g$ of $\varphi_2$. If $F(x)$ is, for example, ``$x = \varphi_0$'', then $g = \varphi_1$. Accordingly, assume $f = \varphi_2$ in {\sc F9}$_{\mathbb{E}}$. So, consider, e.g., $x = \varphi_0$. Such a value for $x$ is different of $\varphi_2$. Besides, $F(\varphi_0)$. That implies $g(\varphi_0) = \varphi_2(\varphi_0) = \varphi_0$. For all remaining values $x\neq\varphi_0$, we know $g$ does not act on $x$. That means $g$ acts solely on $\varphi_0$. And according to {\sc F7}, that function is supposed to be $\varphi_1$. Observe $\varphi_1$ is allowed to have a free occurrence in $F(x)$, according to our Restriction Axiom. Nevertheless $g$ does not act on $\varphi_1$ in our first example. That means either $g(\varphi_1) = \underline{\mathfrak{0}}$ or $g = \varphi_1$. In this case, we have $g = \varphi_1$. Later on we define a membership relationship $\in$ (Definition \ref{pertencercomoagir}) where $x\in f$ iff $f[x]$ and some conditions are imposed over $f$. That entails we can guarantee that in a translation of ZF's Separation Scheme into Flow's language, any free occurrence of $g$ in $F(x)$ will have no impact (in a precise sense). After all, in this first example $F(x)$ is $x = \varphi_0$, while $g$ is $\varphi_1$. For more details see Section \ref{standard}.

Resuming the discussion about restrictions of $\varphi_2$, if $F(x)$ is the formula ``$x = \varphi_0 \vee x = \varphi_1$'', then $g = \varphi_2$. If $F(x)$ is ``$x = x$'', then again $g = \varphi_2$. If $F(x)$ is ``$x \neq x$'', then $g = \varphi_0$. The novelty here, however, happens with the formula $F(x)$ given by ``$x = \varphi_1$''. In that case we have a proper restriction $\gamma$ such that $\gamma\neq \varphi_1$, $\gamma(\gamma) = \gamma$, $\gamma(\varphi_1) = \varphi_1$, and $\gamma(r) = \underline{\mathfrak{0}}$ for any remaining $r$ different of $\varphi_1$ and $\gamma$ itself. Thus, $\gamma$ is a {\em new function\/} whose existence is granted thanks to {\sc F9}$_{\mathbb{E}}$ and no other previous postulate. Besides, {\sc F8} grants $\gamma$ is emergent. Hence, $\sigma_{\gamma} \neq \underline{\mathfrak{0}}$.

\begin{definicao}\label{definindowp}

$z$ is the {\em restricted power of\/} $f\neq\underline{\mathfrak{0}}$ iff $\forall x ((x\neq z \wedge x\neq\underline{\mathfrak{0}})\Rightarrow((z(x) = x \Leftrightarrow x \subseteq f)\wedge (z(x) = \underline{\mathfrak{0}} \Leftrightarrow x\not\subseteq f)))$. We denote $z$ as $\wp(f)$. We adopt the convention $\wp(\underline{\mathfrak{0}}) = \varphi_0$.

\end{definicao}

For example, $\wp(\varphi_0) = \varphi_1$, $\wp(\varphi_1) = \varphi_2$, and $\wp(\varphi_2) = f$, where $f(\varphi_0) = \varphi_0$, $f(\varphi_1) = \varphi_1$, $f(\varphi_2) = \varphi_2$, $f(\gamma) = \gamma$, $f(f) = f$, and $f(x) = \underline{\mathfrak{0}}$ for the remaining values $x$. Recall $\gamma$ acts only on $\varphi_1$ and $\gamma(\varphi_1) = \varphi_1$.

Observe $z = \wp(f)$ is a restriction of $\underline{\mathfrak{1}}$, even if $f$ is not. Besides, $\wp$ is {\em not\/} a function, but a metalinguistic symbol which helps us to abbreviate the formula $z = \wp(f)$ given by the definition above. Observe as well $\wp(f)\subseteq \mathfrak{p}(f)$, for any $f$. While $\mathfrak{p}(f)$ refers to a function who acts on all terms that lurk $f$, $\wp(f)$ acts on all terms that lurk $f$ as long they are restrictions of $f$.

\begin{observacao}\label{consideracoessobrerestricao}

In a sense, {\sc F9}$_{\mathbb{E}}$ is similar to the Separation Scheme in ZFC, since it states the existence of a unique $g$ obtained from a given $f$ and a formula $F(x)$. Nevertheless, the role of Separation Scheme in ZFC is not limited to grant the existence of subsets of a given set. Thanks to that postulate, ZFC avoids antinomies like Russell's paradox. In our case those antinomies are avoided by means of the simple use of Self-Reference (Observation \ref{Russell}). That is one of the reasons why we do not prohibit free occurrences of $g$ in $F(x)$ (like what happens in ZFC). Actually, if we demanded no free occurrences of $g$ in $F(x)$, we would be unable to obtain some useful restrictions, as we can see in the examples below. Nevertheless, we demand $\forall x ((x\neq g \wedge x\neq f) \Rightarrow ((f(x) = g(x)\wedge F(x))\vee (g(x) = \underline{\mathfrak{0}}\wedge\neg F(x)))$ (Definition \ref{delimitadorarestricao}), which is a weaker condition than the prohibition of occurrences of $g$ in $F(x)$. If we recall Observation \ref{Russell}, we can easily see that, for any $y$ and any $r$, neither $y\big|_{x(x)\neq y}$ nor $y\big|_{x(x)\neq r}$ allow us to get any contradiction in the style of Russell's paradox. Finally, if anyone tries to ``define'' a function $f$ from a formula $F(x)$ without using {\sc F9}$_{\mathbb{E}}$, it is perfectly possible to get a contradiction. For example, we can ``define'' a function $f$ as it follows: $\forall x(f(x) = \varphi_0)$. Since no $f$ acts on $\underline{\mathfrak{0}}$, we have a contradiction, namely, $f(\underline{\mathfrak{0}}) = \varphi_0$. Nevertheless, that would not be a definition at all, but simply a new postulate which is inconsistent with our axioms. Definitions are supposed to be conservative, in the sense of not allowing new theorems \cite{daCosta-01}.

\end{observacao}

\begin{definicao}\label{definesim}

$f\sim g$ iff $\forall t ((f[t]\Leftrightarrow g[t])\wedge ((t\neq f \wedge t \neq g)\Rightarrow f(t) = g(t)))$.

\end{definicao}

If $f \sim g$ but $f\neq g$, we say $f$ and $g$ are {\em clones\/}, as discussed in Observation \ref{observaunicidade}. Regarding their images, clones $f$ and $g$ differ solely on $f$ and $g$: $f(g) = \underline{\mathfrak{0}}$, while $g(g) = g$, and $g(f) = \underline{\mathfrak{0}}$, while $f(f) = f$. It is easy to see that $\sim$ is reflexive.

\begin{teorema}\label{marcionoslibertou}

Let $f$ be emergent. Then, for any $h\neq\underline{\mathfrak{0}}$, $h \sim f$ entails $h = f$.

\end{teorema}

\begin{description}
\item[\sc Proof] From Theorem \ref{trioderestricoes}, $f = f\big|_{t = t}$. Since $h \sim f$, $h = f\big|_{t = t}$ as well. But {\sc F9}$_{\mathbb{E}}$ says any restriction is unique. Thus, $h = f$.
\end{description}

This last theorem grants $\underline{\mathfrak{0}}$ and $\varphi_0$ are the only clones in Flow, at least among emergent functions.

\begin{teorema}
$\forall f (\varphi_0\subseteq f).$
\end{teorema}

\begin{description}
\item[\sc Proof] Theorem \ref{noaction} states $\varphi_0$ is the only function different of $\underline{\mathfrak{0}}$ who do not act on any term. Thus, Definition \ref{restricaodefuncao} grants formula above is proven by vacuity.
\end{description}

\begin{teorema}

For any $x$ and $y$, $x\subset y$ entails $y\not\subset x$.

\end{teorema}

\begin{description}
\item[\sc Proof] If $x$ and $y$ are, respectively, $\varphi_0$ and $\underline{\mathfrak{0}}$, the proof is straightforward, since $\underline{\mathfrak{0}}$ is never any restriction, according to {\sc F9}$_{\mathbb{E}}$ and Definition \ref{restricaodefuncao}. For the remaining cases, $x\subset y$ entails $x\neq y$ and, for any $t$, $x[t]\Rightarrow (y[t]\wedge x(t) = y(t))$. So, there is $t'$ where $y[t']\wedge\neg x[t']$ or $y[t']\wedge x[t']\wedge x(t')\neq y(t')$. In anyone of those cases we have $y\not\subset x$.
\end{description}

\begin{teorema}\label{onlyyou}

$\sigma$ and $\underline{\mathfrak{0}}$ are the only functions $f$ such that $\sigma_f = f$.

\end{teorema}

\begin{description}
\item[\sc Proof] We already know $\sigma_\sigma = \sigma$ ({\sc F2}) and $\sigma_{\underline{\mathfrak{0}}} = \underline{\mathfrak{0}}$ (Theorem \ref{fechamentorenato}). If $f\neq\underline{\mathfrak{0}}$ and $\sigma_f = \underline{\mathfrak{0}}$, then $\sigma_f\neq f$. If $f\neq\sigma$, $f\neq\underline{\mathfrak{0}}$ and $\sigma_f\neq\underline{\mathfrak{0}}$, then {\sc F6} and Definition \ref{defineSigmasigma} demand $\sigma_f\neq f$. Observe as well there is no clone of $\sigma$, according to {\sc F6}, although $\sigma$ is not emergent.
\end{description}

\begin{teorema}

$\forall f(\mathbb{E}(f)\Rightarrow f = \sigma_f\big|_{x\neq f})$.

\end{teorema}

\begin{description}
\item[\sc Proof] $\mathbb{E}(f)$ entails $\sigma_f\neq\underline{\mathfrak{0}}$. So, $\sigma_f$ acts on all terms where $f$ acts. But $\sigma_f$ acts on just one more term: $f$ itself. Since formula $F(x)$ in the restriction above is ``$x\neq f$'' (observe we follow all demands for $F(x)$ in {\sc F9}$_{\mathbb{E}}$), then $\sigma_f\big|_{x\neq f\wedge}$ acts exactly on all terms where $f$ acts. Besides, they share all their images, according to {\sc F9}$_{\mathbb{E}}$. So, $\sigma_f\big|_{x\neq f} = f$, since $f[x]\Rightarrow x\neq f$.
\end{description}

\begin{definicao}\label{definindouniaoarbitraria}

Let $f$ be a restriction of $\underline{\mathfrak{1}}$ where $\forall g\forall t((f[g]\wedge g[t])\Rightarrow \mathbb{E}(t))$. The {\em arbitrary union of all terms $g$ where $f$ acts\/} (or {\em arbitrary union of $f$\/}, for short) is defined as $u = \bigcup_{f[g]} g = \underline{\mathfrak{1}}\big|_{\exists g(f[g]\wedge g[t])}$.

\end{definicao}

The idea of arbitrary union is as it follows. If $f$ acts on any $g$ which acts on any $x$, then $u$ acts on that very same $x$; and if no $g$ acts on a given $x$, then $u$ does not act on that $x$. Particularly, we write $u = g\cup h$ for the case where $f$ acts at most on $g$ and $h$.

\begin{definicao}\label{intersecaoarbitraria}

Let $f$ be a restriction of $\underline{\mathfrak{1}}$ where $\forall g\forall t((f[g]\wedge g[t])\Rightarrow \mathbb{E}(t))$. Then the arbitrary intersection of $f$ is given by
$\bigcap_{f[g]}g = \underline{\mathfrak{1}} \big|_{\forall g(f[g] \Rightarrow g[t])}.$

\end{definicao}

Observe the definition above is equivalent to $\bigcap_{f[g]}g = \left( \bigcup_{f[g]} g\right) \big|_{\forall g(f[g] \Rightarrow g[t])}$. If $f$ acts only on $g$ and $h$, the arbitrary intersection may be written as $g\cap h$. In particular, for any $m$ and $n$ from language $\mathcal L$, $\varphi_m\cap\varphi_n = \varphi_m\circ\varphi_n$.

\begin{teorema}

$\bigcup_{\underline{\mathfrak{0}}[g]} = \bigcup_{\varphi_0[g]} = \bigcup_{\varphi_1[g]} = \varphi_0$.

\end{teorema}

Proof is straightforward. Last identity from last theorem illustrates our previous claim that it is possible the union $u$ be one of the terms where $f$ does act: $\varphi_1$ acts on $\varphi_0$, and $\bigcup_{\varphi_1[g]} = \varphi_0$.

\subsection{Ordinals}

\begin{definicao}\label{novadefinicaodeordinais}

Let $F(t)$ be the formula $t\subseteq\underline{\mathfrak{1}}\wedge \mathbb{E}(t)\wedge \forall r(t[r]\Rightarrow r\subseteq t)\wedge \forall r\forall s((t[r]\wedge t[s])\Rightarrow (r = s \vee r[s] \vee s[r]))\wedge \forall r((r\subseteq t\wedge \exists s(r[s]))\Rightarrow \exists m(r[m]\wedge \forall x(r[x]\Rightarrow\neg m[x])))$. Then, $\varpi = \underline{\mathfrak{1}}\big|_{F(t)}$ is called the {\em ordinal function\/}. Moreover, if $\varpi[t]$, we say $t$ is an {\em ordinal\/}.

\end{definicao}

The ordinal function $\varpi$ is defined from a conjunction of five formulas: any ordinal is supposed to be a restriction of $\underline{\mathfrak{1}}$ {\em and\/} emergent {\em and\/} transitive {\em and\/} totaly ordered by actions {\em and\/} well-ordered by actions (term $m$ above is the {\em least\/} element with respect to actions). The first two formulas are obviously satisfied by $\varphi_0$. The remaining ones are satisfied by vacuity. Thus, $\varphi_0$ is an ordinal. That fact entails any $\varphi_n$ is an ordinal. In order to grant we are talking about ordinals in an usual sense, we prove the next theorems.

\begin{teorema}\label{ordinaisatuamsobreordinais}
Every ordinal acts only on ordinals.
\end{teorema}

\begin{description}
\item[\sc Proof] All we have to prove is $\forall r(\forall t(\varpi[t]\wedge t[r])\Rightarrow \varpi[r])$. Formula $F(t)$ from Definition \ref{novadefinicaodeordinais} is a conjunction of five formulas. So, this proof is divided into five parts, where we assume $t$ is an ordinal: {\sc (i)} we know $(\varpi[t]\wedge t[r])\Rightarrow r\subseteq t$. Since $t\subseteq \underline{\mathfrak{1}}$, then $r\subseteq\underline{\mathfrak{1}}$. That means the first formula in the five factors conjunction $F(t)$ is satisfied when we replace $t$ by $r$. {\sc (ii)} According to {\sc F8}, any restriction $r$ of an emergent function $t$ is emergent. Thus, since $t[r]\Rightarrow r\subseteq t$, then $\mathbb{E}(r)$. That settles the second condition from $F(t)$. {\sc (iii)} We know $(\varpi[t]\wedge t[r])\Rightarrow r\subseteq t$. Thus, if $r[s]$, then $t[s]$. That entails $s\subseteq t$. Now, suppose $s\not\subseteq r$. Then, $\exists y(s[y]\wedge \neg r[y])$. Since $t[r]\wedge t[y]$, then $r=y\vee y[r]\vee r[y]$. Nevertheless, the possibilities $r[y]$ and $r = y$ cannot take place. So, we must have $y[r]$. Therefore, $r[s]\wedge s[y]\wedge y[r]$, where $t$ acts on all of them: $r$, $s$, and $y$. Now, let $m$ be a restriction of $t$ who acts only on $r$, $s$, and $y$. Therefore, last formula of the five factor conjunction $F(t)$ cannot be satisfied for $m$ (there is no {\em least\/} term with respect to actions of $m$ on $r$, $s$, and $y$). In other words, our hypothesis $s\not\subseteq r$ leads to a contradiction. Thus, $s\subseteq r$. {\sc (iv)} We know $(\varpi[t]\wedge t[r])\Rightarrow r\subseteq t$. Thus, $(r[u]\wedge r[v])\Rightarrow (t[u]\wedge t[v])$. So, from the definition of $\varpi$, $u = v \vee u[v] \vee v[u]$. {\sc (v)} We know $(\varpi[t]\wedge t[r])\Rightarrow r\subseteq t$. If $u\subseteq r$, then $u\subseteq t$. So, if $u$ acts on at least one term, then last condition of $F(t)$ is satisfied when we replace $t$ by $r$.
\end{description}

\begin{teorema}
If $t$ is an ordinal, then $\sigma_t$ is an ordinal.
\end{teorema}

\begin{description}
\item[\sc Proof] Analogously to last theorem, we split this proof into five parts. {\sc (i)} If $t$ is an ordinal, then it is emergent. According to {\sc F8}, the $\mathfrak{F}$-successor of any emergent function is emergent. That entails $\sigma_t\neq\underline{\mathfrak{0}}$. {\sc (ii)} If $t\subseteq\underline{\mathfrak{1}}$, then $\sigma_t\subseteq\underline{\mathfrak{1}}$. {\sc (iii)} Since $\sigma_t\neq\underline{\mathfrak{0}}$, then $\sigma_t[t]$ and $\forall r(\sigma_t[r]\Rightarrow r\subseteq\sigma_t)$. Last two parts are straightforward.
\end{description}

\begin{teorema}\label{funcaoordinalehabrangente}
The ordinal function $\varpi$ is comprehensive.
\end{teorema}

\begin{description}
\item[\sc Proof] Suppose $\sigma_{\varpi}\neq\underline{\mathfrak{0}}$. If $F(t)$ is the formula used in the definition of $\varpi$, then $F(\varpi)$. But according to item ({\sc ii}) of Definition \ref{delimitadorarestricao}, that entails $\underline{\mathfrak{1}}(\sigma_{\varpi}) = \underline{\mathfrak{0}}$. Such a condition is satisfied only for $\sigma_{\varpi} = \underline{\mathfrak{0}}$, a contradiction. Therefore, $\varpi$ cannot be emergent, although it acts only on emergent terms.
\end{description}

\begin{definicao}
An ordinal $t$ is a {\em limit ordinal\/} iff $\nexists r(\sigma_r = t)$.
\end{definicao}

In particular, $\varphi_0$ is a limit ordinal.

\begin{teorema}\label{ordinaisbemordenadosporacoes}
Any ordinal is well-ordered with respect to actions.
\end{teorema}

\begin{description}
\item[\sc Proof] We adopt the next convention: if $r$ and $s$ are ordinals, $r$ is {\em lesser than\/} $s$ iff $s[r]$. Trichotomy, transitivity, and the existence of a least element (with respect to actions) are straightforward consequences from Definition \ref{novadefinicaodeordinais}.
\end{description}

\begin{definicao}\label{defineordinalimite}
The least limit ordinal who acts on $\varphi_0$ is denoted by $\omega$. We say $r$ is a {\em finite ordinal\/} iff $\omega[r]$.
\end{definicao}

\begin{teorema}
If $r$ is an ordinal, then there is a limit ordinal $s$ such that $s[r]$.
\end{teorema}

\begin{description}
\item[\sc Proof] Analogous to standard literature \cite{Jech-03}.
\end{description}

\subsection{ZF-sets}

It is convenient to introduce the next concepts.

\begin{definicao}\label{definindototo}

If $\mathbb{E}(f)$ and $f$ acts solely on emergent functions, $Dom_f^{\mathfrak{F}} = \underline{\mathfrak{1}}\big|_{f[t]}$ and $Im_f^{\mathfrak{F}} = \underline{\mathfrak{1}}\big|_{\exists x(f[x]\wedge f(x) = t)}$ are, respectively, the {\em $\mathfrak{F}$-domain\/} and the {\em $\mathfrak{F}$-image\/} of $f$. We denote $f$ as $f:Dom_f^{\mathfrak{F}}\to Im_f^{\mathfrak{F}}$.

\end{definicao}

When we write $f:x\to y$, that means $x = Dom_f^{\mathfrak{F}}$ and $y = Im_f^{\mathfrak{F}}$, as long $f$ is emergent. We write $Dom^{\mathfrak{F}}_f$ and $Im^{\mathfrak{F}}_f$ instead of $Dom_f$ and $Im_f$ due to the fact we want to compare $\mathfrak{F}$-domain and $\mathfrak{F}$-image with the usual notions of domain and image (range) of a function in a theory like ZFU, for example. Details in Subsection \ref{modelo}.

Now, let $\alpha(x,y)$ be a formula where there is at least one occurrence of $x$, one occurrence of $y$, and all of them are free. Then, next formula is an axiom.

\begin{description}

\item[\sc F10$_{\alpha}$ - Creation] $\forall x\exists! y(\alpha(x,y))\Rightarrow \forall f(\mathbb{E}(f)\Rightarrow \exists! g(\forall x\forall y((\alpha(x,y)\wedge x\neq g)\Rightarrow ((f[x]\Rightarrow g(x) = y)\wedge (\neg f[x]\Rightarrow g(x) = \underline{\mathfrak{0}})))\wedge \sigma(g) \neq \underline{\mathfrak{0}}\wedge \sigma(Im_g^{\mathfrak{F}})\neq \underline{\mathfrak{0}}\wedge (f[g]\Rightarrow \neg \alpha (g,f(g)))))$.

\end{description}

We employed the notation $\sigma(g)$ instead of $\sigma_g$ due to limitations in the \LaTeX $\,$ editor of this text, particularly for the case of $\mathfrak{F}$-successor of the $\mathfrak{F}$-image of $g$.

Last postulate is supposed to define functions $g$ from certain formulas $\alpha$ and any $f$, as long $f$ is emergent. We say $g$ {\em is created by $f$ and $\alpha(x,y)$\/}, and denote this by $g = f\big|^{\alpha(x,y)}$ or, simply, $g = f\big|^{\alpha}$. Besides, $g$ is emergent, as well as its $\mathfrak{F}$-image (and consequently its $\mathfrak{F}$-domain). Last conjunction in {\sc F10}$_{\alpha}$ avoids ambiguities in the calculation of $g$. Postulate {\sc F10}$_{\alpha}$ motivates us to employ an alternative notation for functions. When we see fit to do so, we may eventually write $a\marcio{f}b$ meaning $f(a) = b$. One obvious advantage of this notation $a\marcio{x}b$ is that it can be used to easily write down strings like $a_1\marcio{f} a_2\marcio{f} a_3\marcio{f} \cdot\cdot\cdot$, meaning $f(a_1) = a_2$, $f(a_2) = a_3$, and so on. We adopt the next convention: every time we state a finite sequence like either $a_1\marcio{f} a_2\marcio{f} a_3\marcio{f} \cdot\cdot\cdot \marcio{f} a_n\marcio{f}\underline{\mathfrak{0}}$ or $a_1\marcio{f} a_2\marcio{f} a_3\marcio{f} \cdot\cdot\cdot \marcio{f} a_n\marcio{f}a_i$ ($i$ is any value among terms from $1$ to $n$), we assume $f(f) = f$ and $f(r) = \underline{\mathfrak{0}}$ if $r$ is different of all previous terms already stated.

{\sc F10}$_{\alpha}$ allows us to get restrictions by other means besides {\sc F9}$_{\mathbb{E}}$. For example, from {\sc F10}$_{\alpha}$ there is a $g$ such that $\varphi_1\marcio{g} \varphi_2\marcio{g} \varphi_3\marcio{g} \underline{\mathfrak{0}}$. In that case $g\subset\sigma$. Function $g$ could be created from, e.g., $\varphi_4$ and an appropriate formula $\alpha$.

With the aid of {\sc F10}$_{\alpha}$ it is quite easy to show $\mathfrak{F}$-composition is not commutative. As an example, consider $f$ as a function given by $\varphi_1\marcio{f}\varphi_2\marcio{f}\varphi_2$ and $g$ given by $\varphi_2\marcio{g}\varphi_3\marcio{g}\varphi_3$. In that case $g\circ f$ is given by $\varphi_1\marcio{g\circ f}\varphi_3\marcio{g\circ f}\underline{\mathfrak{0}}$ and $\varphi_2\marcio{g\circ f}\varphi_3\marcio{g\circ f}\underline{\mathfrak{0}}$, while $f\circ g$ is simply $\varphi_0$ (which is obviously different of $g\circ f$).

It is worth to remark there are some functions which cannot exist (besides those already discussed), as we can see in the next theorem.

\begin{teorema}\label{proibindoumbarracomt}

There is no $x$ where, for a given $t\neq\underline{\mathfrak{1}}$, $t\marcio{x}\underline{\mathfrak{1}}\marcio{x} t$.

\end{teorema}

\begin{description}
\item[\sc Proof] Suppose there is such a function $x$. Then, from {\sc F5}, $x\circ x$ is given by $t\marcio{x\circ x} t$ and $\underline{\mathfrak{1}}\marcio{x\circ x} \underline{\mathfrak{1}}$. But from Theorem \ref{teoremamarcio2} we have that any function $h$ such that $h(\underline{\mathfrak{1}}) = \underline{\mathfrak{1}}$ entails $h = \underline{\mathfrak{1}}$. Thus, $x\circ x$ is supposed to be $\underline{\mathfrak{1}}$. But since $x$ is different of $\underline{\mathfrak{1}}$, then $x\circ x$ cannot be $\underline{\mathfrak{1}}$, according to Theorem \ref{dezenovedemarco}. That is a contradiction!
\end{description}

\begin{teorema}
$\forall f(\mathbb{E}(f)\Rightarrow \exists! g(g = f\big|_F))$.
\end{teorema}

\begin{description}
\item[\sc Proof] Since $f$ is emergent, we can use {\sc F10}$_{\alpha}$. Let $\alpha(x,y)$ be the next formula: $(F(x)\Rightarrow f(x) = y)\wedge (\neg F(x)\Rightarrow y = \underline{\mathfrak{0}})$. In that case, $f\big|^{\alpha} = f\big|_F = g$. Concerning the uniqueness of $g$, that is granted by {\sc F10}$_{\alpha}$.
\end{description}

Next we use the axioms of Creation and Restriction to define ZF-sets, who are supposed to play the same role sets play in ZF set theory.

\begin{definicao}\label{definefuncoesdevonNeumann}

Let $\alpha(u,v)$ be a formula recursively defined as it follows:

\begin{description}

\item[\sc i] $\alpha(\varphi_0,\varphi_0)$;

\item[\sc ii] $\forall u\forall v(\varpi[u]\Rightarrow (\alpha(u,v)\Rightarrow \alpha(\sigma_u,\wp(v))))$;

\item[\sc iii] $\forall u((\varpi[u]\wedge \nexists x(\sigma_x = u))\Rightarrow\alpha(u,\underline{\mathfrak{1}}\big|_{\exists r\exists v(u[r]\wedge \alpha(r,v)\wedge v[t])}))$;

\item[\sc iv] $\forall u((\neg\varpi[u])\Leftrightarrow \alpha(u,\underline{\mathfrak{0}}))$.

\item[\sc v] $\forall u\forall v\forall w((\alpha(u,v)\wedge \alpha(u,w))\Rightarrow v=w)$.

\end{description}

If $\eta = r\big|^{\alpha}$, where $r$ is an ordinal, we say $\eta$ is a {\em von Neumann function with rank\/} $r$. For the sake of abbreviation, any von Neumann function with rank $r$ is denoted by $\eta_r$.

\end{definicao}

As an example consider $\eta_{\varphi_4}$. That is a function given by $\varphi_0\marcio{\eta_{\varphi_4}} \varphi_0$, $\varphi_1\marcio{\eta_{\varphi_4}} \varphi_1$, $\varphi_2\marcio{\eta_{\varphi_4}} \varphi_2$, $\varphi_3\marcio{\eta_{\varphi_4}} \wp(\varphi_2)$, where $\wp(\varphi_0) = \varphi_1$, $\wp(\varphi_1) = \varphi_2$, and $\wp(\varphi_2)$ is a restriction of $\underline{\mathfrak{1}}$ who acts only on $\varphi_0$, $\varphi_1$, $\varphi_2$ and $\gamma$, where $\gamma$ acts only on $\varphi_1$.

\begin{teorema}\label{imagensdefuncoesdevonneumann}

For any ordinal $r$ and any $t$ such that $r[t]$, $\eta_r(t)$ is emergent and a restriction of $\underline{\mathfrak{1}}$.

\end{teorema}

\begin{description}
\item[\sc Proof] Observe we are not saying that any $\eta_r$ is a restriction of $\underline{\mathfrak{1}}$, but only their $\mathfrak{F}$-images. We prove this by transfinite induction. Let $\alpha$ be the formula from Definition \ref{definefuncoesdevonNeumann}. The first ordinal $r$ is $\varphi_0$. {\sc F10}$_{\alpha}$ trivially says $\varphi_0\big|^{\alpha} = \varphi_0$. That entails $\eta_{\varphi_0} = \varphi_0$. Since $\varphi_0$ does not act on any term, then it satisfies the theorem by vacuity. For the sake of argument, next ordinal $\varphi_1$ is the first one who acts on some term. In that case, $\eta_{\varphi_1} = \varphi_1$. The only term where $\varphi_1$ acts is $\varphi_0$. Besides, $\varphi_1(\varphi_0) = \varphi_0$. And $\varphi_0$ is emergent and a restriction of $\underline{\mathfrak{1}}$. Now, suppose all non-$\underline{\mathfrak{0}}$ images of a given $\eta_r$ are emergent and restrictions of $\underline{\mathfrak{1}}$. As we could see, that takes place with the first two ordinals. Since $r$ is an ordinal and any ordinal acts only on ordinals (Theorem \ref{ordinaisatuamsobreordinais}), then, $\eta_r(\sigma_t) = \wp(\eta_r(t))$, for any ordinal $t$ where $r$ acts, as long $r$ acts on $\sigma_t$ besides $t$. From {\sc F8} (see function $h$ in that axiom), $\wp(\eta_r(t))$ is emergent. Besides, from the very definition of restricted power (Definition \ref{definindowp}), $\wp(\eta_r(t))$ is a restriction of $\underline{\mathfrak{1}}$. An analogous result holds for $\eta_{\sigma_r}$, since $\sigma_r$ acts on all terms where $r$ acts and on $r$ itself. In the case where $u$ is a limit ordinal (item {\sc iii} of Definition \ref{definefuncoesdevonNeumann}), then $\eta_{\sigma_u}(u) = \underline{\mathfrak{1}}\big|_{\exists r\exists v(u[r]\wedge \alpha(r,v)\wedge v[t])}$. But, according to {\sc F8} (see function $i$ in that postulate), that last term is emergent. Besides, it is trivially a restriction of $\underline{\mathfrak{1}}$.
\end{description}

\begin{definicao}\label{defineZFconjunto}
{\sc (i)} $\nu = \underline{\mathfrak{1}}\big|_{\exists r(\varpi[r]\wedge\eta_{\sigma_r}(r) = t)}$ is the {\em von Neumann universe\/}; {\sc (ii)} $\forall u(\mathbb{Z}(u)\Leftrightarrow \exists t(\nu[t]\wedge t[u]))$; we read $\mathbb{Z}(u)$ as `{\em $u$ is a ZF-set\/}'.
\end{definicao}

\begin{teorema}\label{todoordinalehZFconjunto}
Every ordinal is a ZF-set.
\end{teorema}

\begin{description}
\item[\sc Proof] For any ordinal $r$, $\eta_{\sigma_{\sigma_r}}(\sigma_r)$ acts on $r$ (although it may act on other terms as well), according to Definition \ref{definefuncoesdevonNeumann}. So, from Definition \ref{defineZFconjunto}, $\mathbb{Z}(r)$.
\end{description}

\begin{teorema}\label{voneumannabrange}
If $\nu$ is the von Neumann universe, then $\sigma_{\nu} = \underline{\mathfrak{0}}$.
\end{teorema}

\begin{description}
\item[\sc Proof] From Theorem \ref{todoordinalehZFconjunto}, $\varpi\subset \nu$. From Theorem \ref{funcaoordinalehabrangente}, $\sigma_{\varpi} = \underline{\mathfrak{0}}$. So, from Theorem \ref{extensaodeabrangenteehabrangente}, $\nu$ is comprehensive. Suppose $\sigma_{\nu}\neq\underline{\mathfrak{0}}$. Then, $\nu$ is emergent, since it acts only on terms with non-$\underline{\mathfrak{0}}$ $\mathfrak{F}$-successor. Thus, from {\sc F8}, every restriction of $\nu$ has a non-$\underline{\mathfrak{0}}$ $\mathfrak{F}$-successor. But $\varpi$ is a restriction of $\nu$ and $\sigma_{\varpi} = \underline{\mathfrak{0}}$, a contradiction. Therefore, $\sigma_{\nu} = \underline{\mathfrak{0}}$.
\end{description}

\begin{teorema}\label{legitimarank}
$\forall u(\mathbb{Z}(u)\Leftrightarrow \exists r(\varpi[r]\wedge u\subseteq\eta_{\sigma_r}(r)))$.
\end{teorema}

\begin{description}
\item[\sc Proof] From Definition \ref{defineZFconjunto}, $\mathbb{Z}(u)\Leftrightarrow \exists r\exists t(\varpi[r]\wedge \eta_{\sigma_r}(r)[t]\wedge t[u])$. From Definition \ref{definefuncoesdevonNeumann}, item {\sc (ii)}, $\eta_{\sigma_{\sigma_r}}(\sigma_r) = \wp(\eta_{\sigma_r}(r))$. Hence, for any $x$, $\eta_{\sigma_{\sigma_r}}(\sigma_r)[x] \Leftrightarrow x\subseteq \eta_{\sigma_r}(r)$.
\end{description}

\begin{definicao}
The {\em rank\/} of a ZF-set $u$ is the least ordinal $r$ such that $u\subseteq\eta_{\sigma_r}(r)$. We denote this by $r = \rank(u)$.
\end{definicao}

As an example, consider a function $a$ such that $b\marcio{a}b$, where $\varphi_1\marcio{b}\varphi_1$. The rank of $a$ is $\varphi_3$, while the rank of $b$ is $\varphi_2$, and the rank of $\varphi_1$ is $\varphi_1$.

\begin{teorema}\label{teoremasobreacoeserank}
$\forall u(\mathbb{Z}(u)\Rightarrow \forall x(u[x]\Rightarrow \rank(u)[\rank(x)]))$.
\end{teorema}

\begin{description}
\item[\sc Proof] Straightforward from Theorem \ref{legitimarank} and the definition of rank of a ZF-set.
\end{description}

Last theorem, in a sense, resembles the definition of rank of a set in \cite{Monk-80}.

\begin{teorema}\label{teoremaquecapturantigadefinicaodezfconjunto}
$\forall f (\mathbb{Z}(f)\Leftrightarrow (\sigma_f\neq\underline{\mathfrak{0}}\wedge f\subseteq\underline{\mathfrak{1}}\wedge\forall x(f[x]\Rightarrow \mathbb{Z}(x)))$.
\end{teorema}

\begin{description}
\item[\sc Proof] The $\Rightarrow$ part is a direct consequence from Theorem \ref{imagensdefuncoesdevonneumann} and Definition \ref{defineZFconjunto}. For the $\Leftarrow$ part, we already know  $f$ is emergent, since $\sigma_f\neq\underline{\mathfrak{0}}$ and $f$ acts only on ZF-sets. That entails we can use {\sc F10}$_\alpha$ to create a function $g$ from $f$, if an adequate formula $\beta$ is used. For that purpose, consider formula $\beta(x,y)$ given by ``$(f[x]\Rightarrow y = \rank (x)) \wedge (\neg f[x]\Rightarrow y = \underline{\mathfrak{0}})$''. According to {\sc F10}$_\alpha$, $\sigma(Im_g^{\mathfrak{F}})\neq \underline{\mathfrak{0}}$. On the other hand, $Im_g^{\mathfrak{F}}$ is a restriction of $\varpi$. Therefore, $Im_g^{\mathfrak{F}}$ is a ZF-set (Theorem \ref{todoordinalehZFconjunto}). Now, let $r = \rank (Im_g^{\mathfrak{F}})$. That entails $r$ acts on the rank of any term $x$ where $f$ acts, i.e., $\forall x(f[x]\Rightarrow r[\rank (x)])$. Now, suppose $f$ is no restriction of $\eta_{\sigma_r}(r)$ (see Theorem \ref{legitimarank}). That entails there is $u$ such that $f[u]\wedge \neg \eta_{\sigma_r}(r)[u]$ (remember $u$ is supposed to be a ZF-set). Since $\neg \eta_{\sigma_r}(r)[u]$, then $\rank (u)[r]$, which entails $\neg r[\rank (u)]$. That contradicts the fact that $\forall x(f[x]\Rightarrow r[\rank (x)])$. Thus, indeed $f\subseteq\eta_{\sigma_r}(r)$. From Theorem \ref{legitimarank}, $f$ is a ZF-set.
\end{description}

\begin{teorema}\label{classedetodosZF}
$u = \underline{\mathfrak{1}}\big|_{\mathbb{Z}(t)}$ is comprehensive.
\end{teorema}

\begin{description}
\item[\sc Proof] Suppose $\sigma_u\neq\underline{\mathfrak{0}}$. Then, according to Theorem \ref{teoremaquecapturantigadefinicaodezfconjunto}, $u$ is a ZF-set. Thus, $\sigma_u = u$, a contradiction.
\end{description}

\begin{definicao}\label{orderedpair}
$f$ is an {\em ordered pair\/} $(a,b)$, with both values $a$ and $b$ different of $\underline{\mathfrak{0}}$, iff there are $\alpha$ and $\beta$ such that $\alpha\neq f$, $\beta\neq f$, $\alpha\neq a$, $\beta\neq b$ and
$$f(x) = \left\{ \begin{array}{cl}
\alpha & \mbox{if}\; x = \alpha\\
\beta & \mbox{if}\; x = \beta\\
\underline{\mathfrak{0}} & \mbox{if}\; x\neq f \wedge x\neq\alpha \wedge x\neq \beta
\end{array}
\right.$$
\noindent
where $\alpha(a) = a$, $\alpha(x) = \underline{\mathfrak{0}}$ if $x$ is neither $a$ nor $\alpha$, $\beta(a) = a$, $\beta(b) = b$, $\beta(x) = \underline{\mathfrak{0}}$ if $x$ is neither $a$ nor $b$ or $\beta$.

\end{definicao}

We have two kinds of ordered pairs: those where $\alpha\neq b$ (first kind) and those where $\alpha = b$ (second kind). Their respective diagrams are as it follows:

\begin{picture}(60,65)

\put(3,25){\framebox(100,30)}


\put(10,45){$f$}

\put(10,11){\framebox(40,30)}

\put(40,43){$\boldmath \curvearrowleft$}

\put(15,32){$\alpha$}

\put(55,11){\framebox(40,30)}

\put(85,43){$\boldmath \curvearrowleft$}

\put(60,31){$\beta$}

\put(60,12){$a$}

\put(58,16){$\boldmath \curvearrowleft$}

\put(13,16){$\boldmath \curvearrowleft$}

\put(77,17){$\boldmath \curvearrowleft$}

\put(80,12){$b$}

\put(15,12){$a$}

\put(203,23){\framebox(100,40)}


\put(210,53){$g$}

\put(220,5){\framebox(75,45)}

\put(240,8){\framebox(50,30)}

\put(280,40){$\boldmath \curvearrowleft$}

\put(285,51){$\boldmath \curvearrowleft$}

\put(225,39){$\beta$}

\put(244,28){$\alpha = b$}

\put(280,10){$a$}

\put(278,15){$\boldmath \curvearrowleft$}

\end{picture}

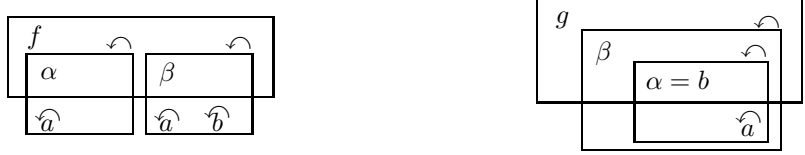
{\sc Figure 2:} Diagrams of $f = (a,b)$ and $g = (a,b)$ of the first and second kind.

Left diagram above concerns the case where $f$ acts only on $\alpha$ and $\beta$, while $\alpha$ acts only on $a$, and $\beta$ acts only on $a$ and $b$. In the particular case where $a = b$, we have $\alpha = \beta$, and the ordered pair $f$ is $(a,a)$. So, $(a,a)$ is a function $f$ which acts solely on $\alpha$, while $\alpha$ acts solely on $a$. Observe that $f(a) = f(b) = \underline{\mathfrak{0}}$ ($f$ never acts neither on $a$ nor on $b$) if $f$ is an ordered pair of the first kind. Thus, $f$ is $(a,b)$ iff $f$ acts only on $\alpha$ and $\beta$, which act, respectively, only on $a$ and only on $a$ and $b$. To get $(b,a)$, all we have to do is to exchange $\alpha$ by a function $\alpha'$ which acts only on $b$. Our definition is obviously inspired on the standard notion due to Kuratowski. In standard set theory $(a,b)$ is a set $\{\{a\}, \{a,b\}\}$ such that neither $a$ nor $b$ belong to $(a,b)$. In Flow, on the other hand, an ordered pair $(a,b)$ of the first kind is a function which does not act neither on $a$ nor on $b$.

Nevertheless, the second kind of ordered pair shows our approach is not equivalent to Kuratowki's. In the case where $\alpha = b$, we have the diagram to the right of {\sc Figure 2}. This non-Kuratowskian ordered pair $g = (a,b)$ acts on $b$, although it does not act on $a$. And no Kuratowskian ordered pair $(a,b)$ ever acts on either $a$ or $b$. So, in the general case, no ordered pair $(a,b)$ ever acts on $a$. That means:

\begin{teorema}\label{aa}

Any ordered pair $(a,a)$ is Kuratowskian.

\end{teorema}

The proof is straightforward. Since any ordered pair $(a,b)$ is a function, for the sake of abbreviation we write $x(a,b)$ for $x((a,b))$, for a given function $x$.

The reader can observe, from {\sc Figure 2}, that in a non-Kuratowskian ordered pair $f$, $\beta$ is the non-$\underline{\mathfrak{0}}$ $\mathfrak{F}$-successor of $\alpha$. In other words:

\begin{teorema}

$f = (a,b)$ is a non-Kuratowskian ordered pair iff $f$ acts only on a function $\alpha$ - which, in its turn, acts on one single term $a$ - and on its $\mathfrak{F}$-successor $\sigma_{\alpha}$, where $\sigma_{\alpha}\neq\underline{\mathfrak{0}}$.

\end{teorema}

\begin{description}
\item[\sc Proof] If $f = (a,b)$, then $f$ acts at most on $\alpha$ and $\beta$, where $\alpha$ acts only on $a$, and $\beta$ acts only on $a$ and $b$. Suppose $f = (a,b)$ is  non-Kuratowskian. Then, $a\neq b$ (Theorem \ref{aa}); and $f = (a,\alpha)$, since $\alpha = b$. But $a\neq b$ entails $a\neq\alpha$. And since $\alpha$ acts on $a$, then $\beta$ acts on two terms: $\alpha$ and $a$. But that is the condition given for granting $\beta = \sigma_{\alpha}$, where $\sigma_{\alpha}\neq\underline{\mathfrak{0}}$. Finally, if $f$ acts only on $\alpha$ and $\sigma_{\alpha}$, where $\alpha$ acts only on $a$ and $\sigma_{\alpha}\neq\underline{\mathfrak{0}}$, then $f = (a,\alpha)$ is non-Kuratowskian.
\end{description}

Suppose $f = \underline{\mathfrak{1}}\big|_{x = \varphi_3 \vee x = \varphi_4}$. In that case $f$ acts on $\varphi_3$ and on its $\mathfrak{F}$-successor $\varphi_4$. Nevertheless, $\varphi_3$ does not act on just one single term. Thus, such an $f$ is not an ordered pair, let alone a non-Kuratowskian ordered pair.

\begin{teorema}

$(a,b) = (c,d)$ iff $a = c$ and $b = d$.

\end{teorema}

\begin{description}
\item[\sc Proof] Straightforward from Definition \ref{orderedpair} and {\sc F9}$_{\mathbb{E}}$.
\end{description}

\begin{teorema}\label{teoreminhaparordenado}
If $a$ and $b$ are both emergent, then there is $f = (a,b)$.
\end{teorema}

\begin{description}
\item[\sc Proof] From {\sc F9}$_{\mathbb{E}}$ we define the proper restriction $\beta$ of $\underline{\mathfrak{1}}$ for ``$t = a \vee t = b$'' as formula $F(t)$. So, $\beta (a) = a$, $\beta (b) = b$, $\beta (\beta) = \beta$, and $\beta (x) = \underline{\mathfrak{0}}$ for all remaining values of $x$. Analogously, the proper restriction $\alpha$ of $\underline{\mathfrak{1}}$ for ``$t = a$''gives us $\alpha (a) = a$, $\alpha (\alpha) = \alpha$, and $\alpha (x) = \underline{\mathfrak{0}}$ for the remaining values of $x$. And, from {\sc F9}$_{\mathbb{E}}$, we have $\alpha\neq a$. Analogously, we have $\beta\neq a$ and $\beta\neq b$. Finally, the proper restriction $f$ of $\underline{\mathfrak{1}}$ for ``$t = \alpha \vee t = \beta$ gives us $f(\alpha) = \alpha$, $f(\beta) = \beta$, $f(f) = f$, and $f(x) = \underline{\mathfrak{0}}$ for all the remaining values of $x$. Besides, $f\neq \alpha$ and $f\neq \beta$. But function $f$ is exactly that one in Definition \ref{orderedpair}. Hence, $f = (a,b)$.
\end{description}

To get an ordered pair $(\underline{\mathfrak{0}},\underline{\mathfrak{0}})$, all we have to do is to consider $\varphi_1$, which acts only on $\varphi_0$. So, $\varphi_1=_{def}(\underline{\mathfrak{0}},\underline{\mathfrak{0}})$. If $f$ acts only on $\alpha$ and $\varphi_0$, where $\alpha$ acts only on $a\neq\underline{\mathfrak{0}}$, then we adopt the convention $f = (a,\underline{\mathfrak{0}})$. For example, $\varphi_2$ is the ordered pair $(\varphi_0,\underline{\mathfrak{0}})$. There are no ordered pairs of the form $(\underline{\mathfrak{0}},b)$, where $b\neq\underline{\mathfrak{0}}$.

\begin{teorema}\label{potenciaZF}

If $f$ is a ZF-set, then $\wp(f)$ is a ZF-set.

\end{teorema}

\begin{description}
\item[\sc Proof] If $f$ is a ZF-set, then it is emergent. So, $f$ acts only on emergent functions and {\sc F8} says any $t$ who lurks $f$ is emergent. Now, let $h = \underline{\mathfrak{1}}\big|_{t\trianglelefteq f}$. Then, $h[t]\Leftrightarrow t\trianglelefteq f$. So, $h = \mathfrak{p}(f)$ (Definition \ref{definindomaximapotencia}). But {\sc F8} also says such an $h$ is emergent. From Definition \ref{definindowp}, $p = \wp(f) = \underline{\mathfrak{1}}\big|_{t\subseteq f}$. But $p$ lurks $h$. So, once again from {\sc F8}, $\sigma_p\neq\underline{\mathfrak{0}}$.
\end{description}

\begin{definicao}\label{localmenteinjetiva}

$f$ is {\em injective\/} iff $\forall r \forall s ((f[r]\wedge f[s] \wedge r\neq s )\Rightarrow f(r)\neq f(s))$. We denote this by $\mathbb{I}(f)$.

\end{definicao}

For example, any $\varphi_n$ is injective.

\begin{definicao}\label{pertencercomoagir}
$x\in_{\mathbb{Z}} f$ iff $f\subset\underline{\mathfrak{1}}\wedge \mathbb{Z}(x)\wedge f[x]$.
\end{definicao}

The negation of formula $x\in_{\mathbb{Z}} f$ is abbreviated as $x\not\in_{\mathbb{Z}} f$. We read $x\in_{\mathbb{Z}} f$ as `{\em $x$ $\mathbb{Z}$-belongs to $f$\/}' or `{\em $x$ is a $\mathbb{Z}$-member of $f$\/}'. The symbol $\in_{\mathbb{Z}}$ is called {\em $\mathbb{Z}$-membership relation\/}.

\begin{teorema}
{\sc (i)} $\forall x (x\not\in_{\mathbb{Z}} x)$; {\sc (ii)} $\forall x(\underline{\mathfrak{0}}\not\in_{\mathbb{Z}} x)$; {\sc (iii)} $\forall x (\underline{\mathfrak{1}}\not\in_{\mathbb{Z}} x)$.
\end{teorema}

\begin{description}
\item[\sc Proof] Item {\sc (i)} is consequence from {\sc F2}, since no $x$ acts on itself. Item {\sc (ii)} is consequence from Theorem \ref{fechamentorenato}: no $x$ acts on $\underline{\mathfrak{0}}$. Item {\sc (iii)} is immediate, since $\underline{\mathfrak{1}}$ is not a ZF-set.
\end{description}

At first glance, item {\sc (i)} looks like an evidence that Flow is well-founded. Nevertheless, the issue of regularity is a little more subtle than that, as we can see later on.

$\{f,g\}=_{def}\underline{\mathfrak{1}}\big|_{x=f\vee x = g}$, where $\mathbb{E}(f) \wedge \mathbb{E}(g)$. So, {\em a pair\/} $\{f,g\}$ is a restriction of $\underline{\mathfrak{1}}$ which acts only on $f$ and $g$, as long they are both emergent.

\begin{teorema}\label{provandopares}

Let $\mathbb{E}(f) \wedge \mathbb{E}(g)$. If $u = \{ f,g\}$ is a pair, then $\mathbb{E}(u)$.

\end{teorema}

\begin{description}
\item[\sc Proof] Let $h = \varphi_2\big|^{\alpha(x,y)}$ be defined from {\sc F10}$_{\alpha}$ by formula $\alpha(x,y)$ given by $(x = \varphi_0\Rightarrow y = f)\wedge (x = \varphi_1 \Rightarrow y = g)\wedge ((x\neq \varphi_0 \wedge x\neq \varphi_1)\Rightarrow y = \underline{\mathfrak{0}})$. According to {\sc F10}$_{\alpha}$, $\sigma_h \neq \underline{\mathfrak{0}}$. Since $h$ acts on emergent functions, then $h$ is emergent. Now, let $u = \underline{\mathfrak{1}}\big|_{t = g(\varphi_0)\vee t = g(\varphi_1)}$. That entails $u$ lurks $h$. And axiom {\sc F8} grants $\sigma_u \neq \underline{\mathfrak{0}}$. Hence, $u$ is emergent.
\end{description}

In particular, any pair of ZF-sets is a ZF-set.

On the left side of the image below we find a Venn diagram of a set $f$ (in the sense of ZF), with its elements $g$, $h$, and $i$. On the right side there is a diagram of a ZF-set $f$ (in the sense of Flow) which acts on $g$, $h$, and $i$. One of the main differences between both pictorial representations is the absence of arrows in the Venn diagram. Venn diagrams do not emphasize the full role of the Extensionality Axiom from ZF: a set is defined either by its elements or by those terms who {\em do not belong\/} to it as well. Within Flow, however, it is explicitly emphasized that anything which is not inside the rectangle of a diagram has an image $\underline{\mathfrak{0}}$. Thus, in a sense, Flow diagrams were always somehow implicit within Venn diagrams. At the end, standard extensional set theories, like ZFC, are simply particular cases of a general theory of functions, as we can formally check in Subsection \ref{venn}.

\begin{picture}(60,45)

\put(33,20){\circle{200}}

\put(13,35){$f$}

\put(23,22){$g$}

\put(40,22){$h$}

\put(32,7){$i$}

\put(153,5){\framebox(100,30)}

\put(160,25){$f$}

\put(172,10){$g$}

\put(170,16){$\boldmath \curvearrowleft$}

\put(200,16){$\boldmath \curvearrowleft$}

\put(230,16){$\boldmath \curvearrowleft$}

\put(202,10){$h$}

\put(232,10){$i$}

\end{picture}

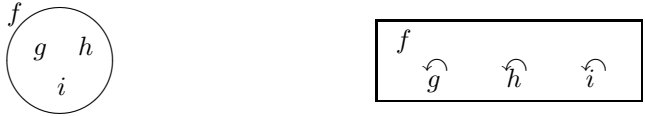
{\sc Figure 3:} Comparison between Venn diagrams and Flow diagrams.

\begin{teorema}

Any inductive function which acts only on ZF-sets is a ZF-set.

\end{teorema}

\begin{description}
\item[\sc Proof] Immediate from the definitions of inductive function and ZF-set.
\end{description}

\begin{definicao}

$f$ is a {\em proper class\/} iff $\forall x(f[x]\Rightarrow \mathbb{Z}(x))\wedge \mathbb{C}(f)$.

\end{definicao}

In other words, no proper class is a ZF-set.

\begin{teorema}\label{unicidadevazio}
There is one single ZF-set $f$ such that for any $x$, we have $x\not\in_{\mathbb{Z}} f$.
\end{teorema}

\begin{description}
\item[\sc Proof] $f = \varphi_0$. From {\sc F6}, $\varphi_0$ is unique. From Theorem \ref{noaction}, $\varphi_0$ is the only term different of $\underline{\mathfrak{0}}$ (recall $\underline{\mathfrak{0}}$ is not a ZF-set) which satisfies the condition given above. In other words, $\varphi_0$ is the {\em empty\/} ZF-set, which can be denoted by $\emptyset$.
\end{description}

\subsection{Atoms}

Next we introduce a generalization of ZF-sets which is quite helpful for us.

\begin{definicao}
\begin{description}
\item[\sc i] $\mathbb{Z}_{\varphi_0}(f)$ iff $\mathbb{Z}(f)$;

\item[\sc ii] $\mathbb{Z}_{\varphi_1}(f)$ iff $f\not\subseteq\underline{\mathfrak{1}}\wedge \mathbb{E}(f)\wedge \forall x(f[x]\Rightarrow (\mathbb{Z}_{\varphi_0}(x)\wedge\mathbb{Z}_{\varphi_0}(f(x))))$;

\item[\sc iii] $\mathbb{Z}_{\sigma_r}(f)$ iff $\varpi[r]\wedge r[\varphi_0]\wedge\mathbb{E}(f) \wedge \exists x (f[x]\wedge (\mathbb{Z}_r(x)\vee \mathbb{Z}_r(f(x))))\wedge \forall y((f[y]\wedge \neg \mathbb{Z}_r(y))\Rightarrow \exists s \exists t(r[s]\wedge r[t]\wedge\mathbb{Z}_s(y)\wedge (\mathbb{Z}_t(f(y))\vee \mathbb{Z}_r(f(y)))))$.

\item[\sc iv] Let $s$ be a limit ordinal different of $\varphi_0$. Then, $\mathbb{Z}_s(f)$ iff $\mathbb{E}(f)\wedge \forall x (f[x]\Rightarrow \exists r\exists t(s[r]\wedge s[t]\wedge \mathbb{Z}_r(x)\wedge \mathbb{Z}_t(f(x))\wedge \exists u\exists y(u[r]\wedge u[t]\wedge s[u]\wedge f[y]\wedge\mathbb{Z}_u(y))))$.
\end{description}
\noindent
If $r$ is an ordinal, we read $\mathbb{Z}_r(f)$ as `$f$ is a {\em ZF-emergent function of degree\/} $r$' or, simply, `{\em ZF-emergent\/}', if the degree $r$ is irrelevant.
\end{definicao}

The only ZF-emergent functions which are required to be no restrictions of $\underline{\mathfrak{1}}$ are those with degree $\varphi_1$. Nevertheless, for our purposes in this Subsection, we are mostly interested on ZF-emergent functions which are no restrictions of $\underline{\mathfrak{1}}$, regardless of their degree, according to next definition.

\begin{definicao}\label{definindoatomos}
$\mathbb{A}(f)$ iff $\exists r (\varpi[r]\wedge \mathbb{Z}_r(f) \wedge f\not\subseteq \underline{\mathfrak{1}})$. We read $\mathbb{A}(f)$ as `{\em $f$ is an atom\/}'.
\end{definicao}

In other words, atoms are ZF-emergent functions who are no restriction of $\underline{\mathfrak{1}}$.

\begin{definicao}\label{tijolo}
Let $a\neq\varphi_0$ be an emergent function who acts only on atoms and such that $a\subset\underline{\mathfrak{1}}$. In that case we say $a$ is a {\em brick of atoms\/}, and denote this by $\mathbb{B}(a)$. If we change item {\sc i} of Definition \ref{definefuncoesdevonNeumann} to $\alpha(\varphi_0,a)$, while keeping all remaining items, then $\eta^a_r = r\big|^{\alpha}$ is an {\em atomic von Neumann function\/}, where $r$ is an ordinal. Besides, $\nu^a = \underline{\mathfrak{1}}\big|_{\exists r(\varpi[r]\wedge\eta^a_{\sigma_r}(r) = t)}$ is an {\em atomic von Neumann universe\/}. Furthermore, $\forall u(\mathbb{M}_a(u)\Leftrightarrow \exists t(\nu^a[t]\wedge t[u]))$; we read $\mathbb{M}_a(u)$ as `{\em $u$ is a Menge in $\nu^a$\/}' or simply `{\em $u$ is a Menge\/}' if there is no risk of confusion.
\end{definicao}

Obviously, the name Menge was borrowed from Cantor's {\em Mengenlehre\/}. Observe no Menge is an atom and no atom is a Menge. Besides, function $a$ of last definition works, for all practical purposes, as the set of all atoms within the context of our atomic von Neumann universe.

\begin{teorema}
For any brick of atoms $a$, $\sigma(\nu^a) = \sigma\left(\underline{\mathfrak{1}}\big|_{\mathbb{M}_a(t)}\right) = \underline{\mathfrak{0}}$.
\end{teorema}

\begin{description}
\item[\sc Proof] Analogous to the proofs of Theorems \ref{voneumannabrange} and \ref{classedetodosZF}.
\end{description}

\begin{teorema}\label{todoZehMenge}
$\forall a\forall f((\mathbb{B}(a) \wedge \mathbb{Z}(f))\Rightarrow \mathbb{M}_a(f))$.
\end{teorema}

The proof of last theorem is trivial. It means every ZF-set is a Menge in any atomic von Neumann universe. Obviously there may be a Menge who is not a ZF-set (recall every ZF-set is a ZF-emergent of degree $\varphi_0$).

Our atoms are supposed to mimic Urelemente from a theory in the style of ZFU (ZF with Urelemente), while certain ZF-emergent functions $f$ such that $f\subset \underline{\mathfrak{1}}$ are supposed to play the role of `collections' or `sets'. In order to do that we need to introduce a specialized concept of `membership' which grants all atoms are empty. Thus:

\begin{definicao}\label{pertinenciasemdecisao}
{\sc (i)} $x\in y$ iff $\mathbb{M}_a(y)\wedge y[x]$, if $a$ is a brick of atoms; {\sc (ii)} $x\not\in y$ iff $\neg(x\in y)$. We read $x\in y$ as `{\em $x$ is an element of $y$\/}' or `{\em $x$ belongs to $y$\/}'.
\end{definicao}

In other words, $x\in y$ iff $y$ is a Menge and $x$ is either an atom or a ZF-set such that $y[x]$ or a Menge such that $y[x]$ and so on. After all, any Menge who acts on some term, acts either on an atom or on a ZF-set or on another Menge. Since no atom is a Menge, then every atom $y$ is empty, in the sense that, for any $x$, we have $x\not\in y$. That means our ZF-sets correspond to pure sets.

From an intuitive point of view, although our atoms are empty (with respect to $\in$), they are not ontologically equivalent to the classical notion of Urelemente. Urelemente, from German, are supposed to be {\em primordial elements\/}, i.e., `unbreakable into minor parts'. However, our atoms have some sort of `inner structure', since they are functions who act on other functions. For example, if $x$ is a function such that $\varphi_0\marcio{x}\varphi_1\marcio{x}\varphi_2\marcio{x}\underline{\mathfrak{0}}$, and $y$ is given by $\varphi_2\marcio{y}\varphi_1\marcio{y}\varphi_0\marcio{y}\underline{\mathfrak{0}}$, they are both atoms. But $x\circ y$ is a ZF-set whose only elements are the ZF-sets $\varphi_1$ and $\varphi_2$. No version of ZFU in literature suggests the possibility of compositions among atoms who are able to produce `pure sets'. Besides, our atoms may be `ranked' by means of degrees: the degree of an atom $a$ is the least ordinal $r$ such that $\mathbb{Z}_r(a)$. That is why we prefer to avoid the terms `Urelement' and `Urelemente' for describing empty terms who are not Menge. On the other hand, the term `atom' is not appropriate either, since its etymology is attached to the notion of indivisibility. Nevertheless, for our purposes in this paper, we refer to our empty terms as atoms, as long they are not Menge but can be elements of a Menge.

If $x$ is a ZF-emergent of degree $\varphi_0$, then it is a Menge in any von Neumann universe (Theorem \ref{todoZehMenge}). If $x$ is a ZF-emergent of degree $\varphi_1$, then it is an atom. Therefore, it is an atom in any atomic von Neumann universe $\nu^a$ where $a$ acts on $x$. If $x$ is a ZF-emergent of degree greater than $\varphi_1$, then either $x\not\subset\underline{\mathfrak{1}}$ or $x\subset\underline{\mathfrak{1}}$. In the first case, $x$ is an atom in $\nu^a$, for any $a$ which acts on $x$. Finally, if $x\subset\underline{\mathfrak{1}}$, then $\forall t(x[t]\Rightarrow (t\subset\underline{\mathfrak{1}}\vee t\not\subset\underline{\mathfrak{1}}))$. Thus, either $b = x\big|_{t\not\subseteq\underline{\mathfrak{1}}}$ is a brick of atoms or $b = \varphi_0$. In the first case $b$ is a Menge in $\nu^c$ for any $c$ such that $b\subseteq c$. In the second case Theorem \ref{todoZehMenge} once again grants $\varphi_0$ is a Menge in any atomic von Neumann universe. All of this points to an interesting phenomenon, namely, ZF-emergent functions of degree $\varphi_1$ constitute some sort of frontier between collections with atoms and collections without atoms.

\subsection{Interpreting ZFU}

Our main goal here is to introduce a minimum set of ingredients to build a model of ZF set theory with atoms (ZFU) where we have PP but not AC. It is well known that permutation models of ZFU are easily capable of violating the Axiom of Choice thanks to one single fact: the axioms of ZFU cannot distinguish among atoms, in the sense that there are non-trivial $\in$-automorphisms on the universe of ZFU. Although Definition \ref{definindoatomos} allows the existence of a vast number of atoms, we are interested only on a specific range of them. But first we need a new concept.

\begin{definicao}
Let $x$ be any function. The {\em kernel\/} of $x$ is given by $\mathcal{K}(x) = \underline{\mathfrak{1}}\big|_{\mathbb{Z}(t)\wedge (x[t]\vee \exists r(x(r) = t))}$.
\end{definicao}

Last definition is sound, since every ZF-set is emergent. Thus, {\sc F9}$_F$ is applicable.

\begin{teorema}
Let $\nu^a$ be an atomic von Neumann universe. Then, $\mathcal{K}(\nu^a) = \nu$, where $\nu$ is the von Neumann universe.
\end{teorema}

\begin{description}
\item[\sc Proof] Straightforward from Theorem \ref{todoZehMenge}.
\end{description}

Thus, we have here a generalization of the usual concept of kernel \cite{Jech-03} which is applicable in a nontrivial way to atoms as well, as we can see in the next Observation.

\begin{observacao}\label{exemplosdeatomos}
According to our previous discussions we have a lot of options for choosing a brick $a$ of atoms for defining any atomic von Neumann universe $\nu^a$. That is why we now introduce some simple examples for the sake of illustration. Let $\mu$ be a function obtained from $\varpi\big|_{t = \omega \vee t = \sigma_{\omega}}$ and {\sc F10}$_{\alpha}$ as it follows, where $\omega$ is the limit ordinal who acts only on finite ordinals (Definition \ref{defineordinalimite}): $\omega\marcio{\mu}\sigma_{\omega}\marcio{\mu}\omega$. Let $\mu'$ be obtained from $\varpi\big|_{t = \sigma_{\omega}\vee t = \sigma_{\sigma_{\omega}}}$ and {\sc F10}$_{\alpha}$ as it follows: $\sigma_{\omega}\marcio{\mu'}\sigma_{\sigma_{\omega}}\marcio{\mu'}\sigma_{\omega}$. Obviously, $\mu\neq\mu'$. Let $p$ be a function obtained from $\varphi_3$ and {\sc F10}$_{\alpha}$ as it follows: $\varphi_0\marcio{p}\varphi_1\marcio{p}\varphi_2\marcio{p}\varphi_0$. Using once again {\sc F10}$_{\alpha}$ over $p$ we can obtain the next function $g$ given by $\varphi_1\marcio{g}\varphi_0\marcio{g}\mu\marcio{g}\underline{\mathfrak{0}}$ and $\varphi_2\marcio{g}\mu'\marcio{g}\underline{\mathfrak{0}}$. Analogously there is an $h$ such that $\mu\marcio{h}\varphi_0\marcio{h}\varphi_1\marcio{h}\underline{\mathfrak{0}}$ and $\mu'\marcio{h}\varphi_2\marcio{h}\underline{\mathfrak{0}}$. Hence, there is a function $a_{\varphi_0}$ where $\mu\marcio{a_{\varphi_0}}\varphi_0\marcio{a_{\varphi_0}}\mu'\marcio{a_{\varphi_0}}\mu$. By using analogous arguments we introduce the next atom $\overline{a_{\varphi_0}}$: $\mu\marcio{\overline{a_{\varphi_0}}}\mu'\marcio{\overline{a_{\varphi_0}}}\varphi_0\marcio{\overline{a_{\varphi_0}}}\mu$. Observe $a_{\varphi_0}\neq\overline{a_{\varphi_0}}$. Besides, $a_{\varphi_0}\circ a_{\varphi_0} = \overline{a_{\varphi_0}}$, $\overline{a_{\varphi_0}}\circ\overline{a_{\varphi_0}} = a_{\varphi_0}$, but $\mathcal{K}(a_{\varphi_0}) = \mathcal{K}(\overline{a_{\varphi_0}}) = \mathcal{K}(a_{\varphi_0}\circ a_{\varphi_0}) = \mathcal{K}(\overline{a_{\varphi_0}}\circ \overline{a_{\varphi_0}}) = \mathcal{K}(a_{\varphi_0}\circ\overline{a_{\varphi_0}}) = \mathcal{K}(\overline{a_{\varphi_0}}\circ a_{\varphi_0}) = \varphi_1$. Although both atoms $a_{\varphi_0}$ and $\overline{a_{\varphi_0}}$ are empty (there is no $t$ such that $t$ belongs to either one of them), their kernels are $\varphi_1$, a non-empty ZF-set. We say $a_{\varphi_0}$ and $\overline{a_{\varphi_0}}$ are {\em conjugates\/} of each other, in the sense that they cannot be distinguished from each other neither by means of their kernels nor by means of the kernels of their $\mathfrak{F}$-compositions. So, yes, we are inspired on the original ideas by Abraham A. Fraenkel \cite{Fraenkel-22}, despite the fact that his seminal paper has some well known mistakes \cite{Mostowski-38}. Nevertheless we use an analogous technique to introduce atoms which can be distinguished from $a_{\varphi_0}$ and $\overline{a_{\varphi_0}}$. That is a key point for our main result.
\end{observacao}

\begin{definicao}\label{rankeandocomatomos}
The {\em rank\/} of a Menge $u$ in an atomic von Neumann universe $\nu^a$ is the least ordinal $r$ such that $u\subseteq\eta_{\sigma_r}^a(r)$. We denote this by $r = \rank(u)$. The rank of any atom in an atomic von Neumann universe is $\varphi_0$.
\end{definicao}

Observe this last concept of rank of a Menge in $\nu^a$ is independent from the notion of degree of any atom where $a$ acts. Next we introduce the concept of {\em transitive closure\/}, which is a copy from the usual notion in set theory.

\begin{definicao}
Let $\mathbb{M}_a(x)$, i.e., $x$ is a Menge from an atomic von Neumann universe $\nu^a$. Then: {\sc (i)} $x_{\varphi_0} = x$; {\sc (ii)} $x_{\sigma_r} = \bigcup_{x_r[t]}t$, where $r$ is a finite ordinal; {\sc (iii)} $\mbox{TC}(x) = \bigcup_{\omega[r]}x_r$. We read $\mbox{TC}(x)$ as the `{\em transitive closure\/}' of $x$.
\end{definicao}

\begin{definicao}\label{definindoindiscerniveis}
Let $a$ be a brick of atoms such that $\forall x(a[x]\Rightarrow \exists! y(a[y]\wedge x\neq y\wedge \mathcal{K}(x) = \mathcal{K}(y)\wedge \mathcal{K}(x\circ y) = \mathcal{K}(y\circ x)))$. Then, in the corresponding atomic von Neumann universe $\nu^a$, we have:
\begin{description}
\item[\sc i] If $\mathbb{Z}(x)$ and $y$ is any term, then $x\equiv y$ iff $x = y$;
\item[\sc ii] If $a[x]\wedge a[y]$, then $x\equiv y$ iff $\mathcal{K}(x) = \mathcal{K}(y) \wedge \mathcal{K}(x\circ y) = \mathcal{K}(y\circ x)$;
\item[\sc iii] Let $x = \underline{\mathfrak{1}}\big|_{t=r\vee t=s}$, where $a[r]\wedge a[s] \wedge r\neq s \wedge r\equiv s$; then $\forall y(y\equiv x\Leftrightarrow y = x)$;
\item[\sc iv] Let $f = \{ r\}_x$ be an abbreviation for the formula $x\subseteq a\wedge a[r]\wedge\forall s(f[s]\Leftrightarrow (x[s]\wedge s\equiv r))$; now, let $x\subseteq a \wedge y\subseteq a$, then $x\equiv y$ iff $\forall u(x[u]\Rightarrow \exists! v(v = \{u\}_x\wedge v\equiv \{u\}_y))$;
\item[\sc v] Let $\mathbb{M}_a(x)\wedge\mathbb{M}_a(y)$, then $x\equiv y$ iff $\rank(x) = \rank(y)\wedge \forall r\forall s((\varpi[r]\wedge\varpi[s])\Rightarrow \mbox{TC}(x\big|_{\rank(t) = r \vee \rank(t) = s})\big|_{a[m]} \equiv \mbox{TC}(y\big|_{\rank(t) = r \vee \rank(t) = s})\big|_{a[m]})$, where all occurrences of $m$ in formula $a[m]$ are free;
\item[\sc vi] $\forall x\forall y(x\equiv y \Leftrightarrow y\equiv x)$;
\item[\sc vii] If $a[x]\wedge \mathbb{M}_a(y)$, then $\neg(x\equiv y)$.
\end{description}
We read $x\equiv y$ as `{\em $x$ is indistinguishable of $y$\/}' or `{\em $x$ is indiscernible of $y$\/}'. $x\not\equiv y$ abbreviates $\neg(x\equiv y)$, and we read it as `{\em $x$ is distinguishable of $y$\/}' or `{\em $x$ is discernible of $y$\/}'.
\end{definicao}

Obviously, any $x$ of $\nu^a$ (i.e., either $a[x]$ or $\mathbb{M}_a(x)$) is indistinguishable of itself.

The same criteria used for indiscernible atoms (item {\sc ii}) could be used to ZF-sets. After all, if both $x$ and $y$ are ZF-sets, then $x\circ y = y\circ x$ is simply the intersection between $x$ and $y$. In that case indistinguishability $\equiv$ is simply identity $=$, since $\mathcal{K}(x) = x$, $\mathcal{K}(y) = y$, $\mathcal{K}(x\circ y) = x\circ y$, and $\mathcal{K}(y\circ x) = y\circ x$. Nevertheless, there may exist {\em different\/} atoms which are indistinguishable, like $a_{\varphi_0}$ and $\overline{a_{\varphi_0}}$ in Observation \ref{exemplosdeatomos}. That is a key feature we intend to use in our proof of independence of AC from the Partition Principle within ZFU.

Item {\sc iv} of last definition introduces {\em classes of indiscernible atoms\/}. For instance, let $a_{\varphi_0}$ be the atom $\mu\marcio{a_{\varphi_0}}\varphi_0\marcio{a_{\varphi_0}}\mu'\marcio{a_{\varphi_0}}\mu$. Let $\overline{a_{\varphi_0}}$ be the atom $\mu\marcio{\overline{a_{\varphi_0}}}\mu'\marcio{\overline{a_{\varphi_0}}}\varphi_0\marcio{\overline{a_{\varphi_0}}}\mu$. Those are the same ones from Observation \ref{exemplosdeatomos}. Now, let $a_{\varphi_1}$ be the atom $\mu\marcio{a_{\varphi_1}}\varphi_1\marcio{a_{\varphi_1}}\mu'\marcio{a_{\varphi_1}}\mu$, and $\overline{a_{\varphi_1}}$ be the atom $\mu\marcio{\overline{a_{\varphi_1}}}\mu'\marcio{\overline{a_{\varphi_1}}}\varphi_1\marcio{\overline{a_{\varphi_1}}}\mu$. Observe $a_{\varphi_0}\not\equiv a_{\varphi_1}$. If $x$ is a Menge who acts only on $a_{\varphi_0}$, $\overline{a_{\varphi_0}}$, and $a_{\varphi_1}$, then $\{a_{\varphi_0}\}_x = \{a_{\overline{\varphi_0}}\}_x$, namely, a Menge who acts only on $a_{\varphi_0}$ and $\overline{a_{\varphi_0}}$. On the other hand, $\{a_{\varphi_1}\}_x = \{a_{\overline{\varphi_1}}\}_x$, namely, a Menge who acts solely on $a_{\varphi_1}$. If $a_{\varphi_2}$ is defined in a similar fashion, then $\{a_{\varphi_2}\}_x = \varphi_0$.

Bricks of atoms in the sense of last Definition provide necessary but not sufficient conditions for achieving our goals. That is why we need the next concept. We introduce below our {\em Appropriate Brick of Atoms\/} which is used in our main result of this paper.

\begin{definicao}\label{tijoloapropriado}
For any finite ordinal $r$, let $\alpha_r$ be the atom $r\marcio{\alpha_r}\mu\marcio{\alpha_r}\mu'\marcio{\alpha_r} r$ and $\overline{\alpha_r}$ be the atom $r\marcio{\alpha_r}\mu'\marcio{\alpha_r}\mu\marcio{\alpha_r} r$, where $\mu$ and $\mu'$ are the same functions introduced in Observation \ref{exemplosdeatomos}. Then, $\mathfrak{a} = \underline{\mathfrak{1}}\big|_{\exists r(\omega[r]\Rightarrow (t = \alpha_r\vee t = \overline{\alpha_r}))}$ is the {\em Appropriate Brick of Atoms\/}, where $\omega$ is the limit ordinal who acts only on all finite ordinals (Definition \ref{defineordinalimite}), i.e., all $\varphi_n$. Any $t$ where $\mathfrak{a}$ acts is an {\em appropriate atom\/}. Within this context, $\nu^{\mathfrak{a}}$ is the {\em Appropriate Atomic von Neumann Universe\/}, or {\em Appropriate Universe\/}, for short.
\end{definicao}

That means the Appropriate Brick of Atoms $\mathfrak{a}$ is a Menge in $\nu^{\mathfrak{a}}$. From now on we use the same notation of last definition for appropriate atoms, i.e., anytime we write either $\alpha_r$ or $\overline{\alpha_r}$ we mean both are appropriate atoms. Next theorem is quite important, although its proof is a simple exercise.

\begin{teorema}\label{basedorelativismo}
For any finite ordinals $r$ and $s$, such that $r\neq s$, we have: {\sc (i)} $\mathcal{K}(\alpha_r) = \mathcal{K}(\overline{\alpha_r}) = \mathcal{K}(\alpha_r\circ \alpha_r) = \mathcal{K}(\overline{\alpha_r}\circ \overline{\alpha_r}) = \mathcal{K}(\alpha_r\circ \overline{\alpha_r}) = \mathcal{K}(\overline{\alpha_r}\circ \alpha_r) = \gamma$, where $\gamma$ is given by $r\marcio{\gamma}r$; {\sc (ii)} $\mathcal{K}(\alpha_r\circ \alpha_s) = \mathcal{K}(\overline{\alpha_r}\circ \alpha_s) = \mathcal{K}(\alpha_r\circ \overline{\alpha_s}) = \mathcal{K}(\overline{\alpha_r}\circ \overline{\alpha_s}) = \gamma'$, where $\gamma'$ is given by $r\marcio{\gamma'}r$ and $s\marcio{\gamma'}s$.
\end{teorema}

\begin{description}
\item[\sc Proof] Just a bunch of boring calculations.
\end{description}

\begin{teorema}\label{discerniveisbasicos}
For any finite ordinals $r$ and $s$, such that $r\neq s$, we have: {\sc (i)} $\alpha_r\equiv \overline{\alpha_r}$; {\sc (ii)} $\alpha_r\not\equiv\alpha_s$; {\sc (iii)} $\alpha_r\not\equiv\overline{\alpha_s}$; {\sc (iv)} $\overline{\alpha_r}\not\equiv\alpha_s$; {\sc (v)} $\overline{\alpha_r}\not\equiv\overline{\alpha_s}$.
\end{teorema}

\begin{description}
\item[\sc Proof] For any finite ordinal $r$ we have $\mathcal{K}(\alpha_r) = \mathcal{K}(\overline{\alpha_r}) = \mathcal{K}(\alpha_r\circ \overline{\alpha_r}) = \mathcal{K}(\overline{\alpha_r}\circ \alpha_r) = \gamma$, where $\gamma$ is given by $r\marcio{\gamma}r$. That settles item {\sc i}. For the remaining items observe that $\mathcal{K}(\alpha_s) = \mathcal{K}(\overline{\alpha_s}) = \gamma'$, where $r\marcio{\gamma'}r$ and $s\marcio{\gamma'}s$. Since $r\neq s$, then $\gamma\neq\gamma'$.
\end{description}

\begin{teorema}\label{tambemutilparalema7}
$\forall f(\mathbb{M}_{\mathfrak{a}}(f) \Leftrightarrow (\sigma_f\neq\underline{\mathfrak{0}}\wedge f\subseteq\underline{\mathfrak{1}}\wedge \forall x(f[x]\Rightarrow (\mathfrak{a}[x]\vee\mathbb{M}_{\mathfrak{a}}(x)))))$.
\end{teorema}

\begin{description}
\item[\sc Proof] Analogous to the proof of Theorem \ref{teoremaquecapturantigadefinicaodezfconjunto}, since any Menge in $\nu^{\mathfrak{a}}$ can be ranked.
\end{description}

Next we follow a maneuver somehow inspired on Cantor-Schr\"oder-Bernstein theorem \cite{Kolmogorov-75}.

\begin{definicao}\label{equipotencia}

$g$ and $h$ are {\em equipotent\/}, denoted by $g\cong h$, iff there is $\tau$ where, for any $t$: {\sc (i)} $\tau[t]\Rightarrow\tau(\tau(t)) = t$; {\sc (ii)} $g[t]\Rightarrow (\tau[t] \wedge h[\tau(t)])$; {\sc (iii)} $h[t]\Rightarrow (\tau[t] \wedge g[\tau(t)])$; {\sc (iv)} $\tau[t]\Rightarrow (g[t]\vee h[t])$. We call function $\tau$ a {\em connector of $g/h$\/}.

\end{definicao}

For example, any $\varphi_n$ is equipotent to itself. And a connector of $\varphi_n/\varphi_n$ is $\varphi_n$.

\begin{definicao}\label{funcaodedekindinfinita}

$g$ is {\em Dedekind-finite\/} iff $\exists h(g\cong h \wedge \tau\not\cong g)$, where $\tau$ is a connector of $g/h$ (Definition \ref{equipotencia}); and $g$ is {\em Dedekind-infinite\/} iff it is not Dedekind-finite.

\end{definicao}

As an example, $h = \varphi_5\big|_{t = \varphi_3\vee t = \varphi_4}$ acts only on $\varphi_3$ and $\varphi_4$. So, if we prove the existence of $\tau$ such that $\tau(\varphi_0) = \varphi_3$, $\tau(\varphi_3) = \varphi_0$, $\tau(\varphi_1) = \varphi_4$, and $\tau(\varphi_4) = \varphi_1$, then $h\cong \varphi_2$ with connector $\tau$. Besides, $\tau\not\cong\varphi_2 \wedge \tau\not\cong h$, since $\tau$ acts on four terms, while $\varphi_2$ and $h$ act on two terms each. To prove the existence of such a $\tau$ we use {\sc F10}$_{\alpha}$.

\begin{teorema}

$\cong$ is reflexive among emergent functions.

\end{teorema}

\begin{description}
\item[\sc Proof] If $f$ is $\varphi_0$, assume $\underline{\mathfrak{0}}$ as connector, according to Definition \ref{equipotencia}. For the remaining cases, make $\tau = \underline{\mathfrak{1}}\big|_{f[x]}$ as a connector of $f/f$.
\end{description}

\begin{teorema}\label{equipotenciaquelevaaequipotencia}

$\forall f\forall g ((\mathbb{Z}(f)\wedge \mathbb{Z}(g)\wedge f\cong g)\Rightarrow \sigma_f\cong \sigma_g)$.

\end{teorema}

\begin{description}
\item[\sc Proof] If $f$ and $g$ are ZF-sets, then $f\cup g$ is a ZF-set. That means $\sigma_f\neq\underline{\mathfrak{0}}$, $\sigma_g\neq\underline{\mathfrak{0}}$, and $\sigma_{f\cup g}\neq\underline{\mathfrak{0}}$. If $f\cong g$, then there is a connector $\tau$ of $f/g$. Since $f\cup g$ is uncomprehensive, we can apply axiom {\sc F10}$_{\alpha}$ over $f\cup g$, where formula $\alpha(x,y)$ is $\tau(x) = y$. That entails, from {\sc F10}$_{\alpha}$, that $\sigma_{\tau}\neq\underline{\mathfrak{0}}$. Besides, $\sigma_f$, $\sigma_g$, and $\sigma_{f\cup g}$ are uncomprehensive and ZF-sets. That entails we can define a new connector $\tau'$ for $\sigma_f/\sigma_g$. All we have to do is to apply {\sc F10}$_{\alpha}$ over $\sigma_f\cup\sigma_g$ with formula $\alpha'(x,y)$ given by $\alpha(x,y)$ when $x\neq f$ and $x\neq g$, and such that $\alpha'(f,g)$ and $\alpha'(g,f)$. Hence, $\sigma_f\cong\sigma_g$, where $\tau'$ is the connector of $\sigma_f/\sigma_g$.
\end{description}

\begin{definicao}
$g$ is {\em infinite\/} iff there is no finite ordinal $r$ such that $r\cong g$. Otherwise, we say $g$ is {\em finite\/}.
\end{definicao}

If we recall that $\mathbb{E}$ (emergent), $\mathbb{I}$ (injective), and $\mathbb{B}$ (brick of atoms) are given, respectively, in Definitions \ref{definindoemergente}, \ref{localmenteinjetiva}, and \ref{tijolo}, next formula is our last axiom of this first version of Flow.

\begin{description}

\item[\sc F11 - $\mathfrak{F}$-Partition] $\forall f\forall a((\mathbb{B}(a)\wedge \mathbb{E}(f)\wedge \mathbb{M}_a(Dom_f^{\mathfrak{F}})\wedge \mathbb{M}_a(Im_f^{\mathfrak{F}}))\Rightarrow \exists c(Dom_c^{\mathfrak{F}} = Im_f^{\mathfrak{F}}\wedge Im_c^{\mathfrak{F}}\subseteq Dom_f^{\mathfrak{F}}\wedge \mathbb{I}(c)))$.

\end{description}

\subsection{The model}\label{modelo}

Next we prove that if $a = \mathfrak{a}$, then we are able to introduce a model of ZFU where the Partition Principle holds but the Axiom of Choice fails.

\begin{definicao}
Let $a$ be a brick of atoms. A function $f$ is a {\em ZFU-function\/} in an atomic von Neumann universe $\nu^a$ iff:
\begin{description}
\item[\sc i] $\mathbb{M}_a(f)$;
\item[\sc ii] $\forall t(f[t]\Rightarrow \exists u\exists v((a[u]\vee a[v]\vee \mathbb{M}_a(u) \vee \mathbb{M}_a(v))\wedge t = (u,v)))$;
\item[\sc iii] $\forall u\forall v\forall w((f[(u,v)]\wedge f[(u,w)])\Rightarrow v = w)$.
\end{description}
\end{definicao}

Observe the concept of ZFU-function corresponds to the usual notion of function in a theory like ZFU: a specific set of ordered pairs $(u,v)$ where $u$ and $v$ are either atoms or sets.

\begin{teorema}\label{funcaoZFfuncaocorrespondem}
Let $a$ be a brick of atoms and $f:x\to y$ (Definition \ref{definindototo}) be a function such that both $x$ and $y$ are Mengen of $\nu^a$. Then there is one single ZFU-function $f'$ such that $\forall u\forall v(f[u]\wedge f(u) = v)\Leftrightarrow (f'[(u,v)]\wedge f'((u,v)) = (u,v))$. Conversely, if $f'$ is a ZFU-function in $\nu^a$, then there is a single $f:x\to y$ such that $\mathbb{M}_a(x)\wedge \mathbb{M}_a(y)$ and $\forall u\forall v(f[u]\wedge f(u) = v)\Leftrightarrow (f'[(u,v)]\wedge f'((u,v)) = (u,v))$.
\end{teorema}

\begin{description}
\item[\sc Proof] First part. Let $f:x\to y$ be a function such that $\mathbb{M}_a(x)$ and $\mathbb{M}_a(y)$. Then, $\forall r(f[r]\Rightarrow \exists s(y[s]\wedge f(r) = s))$. Now, let $\alpha(r,t)$ be the next formula: $(f[r]\Rightarrow t = (r,f(r)))\wedge (\neg f[r]\Rightarrow t = \underline{\mathfrak{0}})$. Then, $h = f\big|^{\alpha}$ is a function whose non-$\underline{\mathfrak{0}}$ images $h(r)$ are just ordered pairs $(r,s)$ such that $f(r) = s$. Besides, {\sc F10}$_{\alpha}$ entails $\sigma(Im_{h}^{\mathfrak{F}})\neq\underline{\mathfrak{0}}$. Another way to write this is $f' = Im_{h}^{\mathfrak{F}} = \underline{\mathfrak{1}}\big|_{\forall a\forall b(t = (a,b) \Leftrightarrow (f[a]\wedge f(a) = b))}$. So, $f'$ is a ZFU-function, whose uniqueness is granted by {\sc F9}$_{\mathbb{E}}$. Theorem \ref{tambemutilparalema7} grants $f'$ is a Menge in $\nu^a$, since any ordered pair $(r,s)$ of either atoms or Mengen has a non-$\underline{\mathfrak{0}}$ $\mathfrak{F}$-successor. Second part. Let $f'$ be a ZFU-function in $\nu^a$. Then $\forall u\forall v\forall w((f'[(u,v)]\wedge f'[(u,w)])\Rightarrow v = w)$. Now, let $x = \underline{\mathfrak{1}}\big|_{\forall u(t = u \Leftrightarrow \exists v (f'[(u,v)]))}$. By analogous arguments used in the first part, $\sigma_x\neq\underline{\mathfrak{0}}$. Next, let $\alpha(r,s)$ be the formula $(x[r]\Rightarrow \forall t(f'[(r,t)]\Leftrightarrow t=s))\wedge (\neg x[r]\Rightarrow s = \underline{\mathfrak{0}})$. The definition of ZFU-function grants $\forall r\exists! s(\alpha(r,s))$. That means we can use {\sc F10}$_{\alpha}$. Thus, $f = x\big|^{\alpha}$ is a function $f:x\to y$ such that $\forall r(f[r]\Rightarrow \exists s(y[s]\wedge f(r) = s))$. Besides, $\sigma_y \neq \underline{\mathfrak{0}}$. The uniqueness of $f$ is granted by {\sc F10}$_{\alpha}$. Theorem \ref{tambemutilparalema7} grants $x$ and $y$ are Mengen in $\nu^a$.
\end{description}

Last theorem could be rewritten in terms of an injective function $\mathfrak{f}$ from the `space' of all ZFU-functions of $\nu^a$ to the `space' of all functions $f:x\to y$ where $x$ and $y$ are Mengen of $\nu^a$.

\begin{definicao}\label{correspondencia}
We say $f:x\to y$ is {\em in correspondence\/} to a ZFU-function $f'$ in a given $\nu^a$ iff $\forall u\forall v(f[u]\wedge f(u) = v)\Leftrightarrow (f'[(u,v)]\wedge f'((u,v)) = (u,v))$, where $x$, $y$, and $f'$ are Mengen in $\nu^a$.
\end{definicao}

\begin{definicao}\label{dominioerangedezfu}
Let $f'$ be a ZFU-function. Then, the {\em domain\/} of $f'$ is $Dom_{f'} = \underline{\mathfrak{1}}\big|_{\forall u(t = u \Leftrightarrow \exists v (f'[(u,v)]))}$, and the {\em range\/} of $f'$ is $Im_{f'} = \underline{\mathfrak{1}}\big|_{\forall v(t = v \Leftrightarrow \exists u (f'[(u,v)]))}$.
\end{definicao}

\begin{teorema}\label{compartilhamdominioeimagem}
If $f:x\to y$ is in correspondence to a ZFU-function $f'$, then $Dom^{\mathfrak{F}}_f = Dom_{f'} = x$ and $Im^{\mathfrak{F}}_f = Im_{f'} = y$. Besides, $f\neq f'$.
\end{teorema}

\begin{description}
\item[\sc Proof] See the proof of Theorem \ref{funcaoZFfuncaocorrespondem}.
\end{description}

Next definition provides all necessary and sufficient conditions to introduce the universe of our main model.

\begin{definicao}\label{levy}
Let $F(t)$ be the formula $\exists v(\nu^{\mathfrak{a}}[v]\wedge v[t]\wedge \forall r(\mbox{TC}(t)\big|_{\mathfrak{a}[m]}[r]\Rightarrow \exists s(r\equiv s\wedge r\neq s\wedge \mbox{TC}(t)\big|_{\mathfrak{a}[m]}[s])))$. Then, $\nu^{\mathfrak{a}}_{\dagger} = \nu^{\mathfrak{a}}\big|_{F(t)}$ is the {\em Azriel L\'evy universe\/} or {\em L\'evy universe\/}, for short. If $x$ is a Menge in $\nu^{\mathfrak{a}}_{\dagger}$, we denote this by $\mathbb{M}_{\mathfrak{a}}^{\dagger}(x)$.
\end{definicao}

The L\'evy universe is a restricted appropriate atomic von Neumann universe in the sense that $\nu^{\mathfrak{a}}_{\dagger} \subset \nu^{\mathfrak{a}}$. Last definition grants there is no Menge $f$ in $\nu^{\mathfrak{a}}_{\dagger}$ which acts on a single $u$ if there is any $v$ in $\nu^{\mathfrak{a}}$ where $v\equiv u$ and $v\neq u$. For example, Flow allows us to define a Menge $f$ in $\nu^{\mathfrak{a}}$ such that $\alpha_r\marcio{f}\alpha_r$, where $\alpha_r$ is an appropriate atom. Nevertheless, $\nu^{\mathfrak{a}}_{\dagger}$ does not act on any $t$ such that $t$ acts on that $f$. If $g$ is a Menge in $\nu^{\mathfrak{a}}_{\dagger}$ who acts on any $\alpha_r$, then $g$ acts on $\overline{\alpha_r}$ as well, since $\overline{\alpha_r}$ is indiscernible from $\alpha_r$ (although $\alpha_r\neq\overline{\alpha_r}$).

\begin{teorema}\label{colecoesindiscerniveisaoiguais}
$\forall x\forall y((\mathbb{M}_{\mathfrak{a}}^{\dagger}(x)\wedge \mathbb{M}_{\mathfrak{a}}^{\dagger}(y)\wedge x\equiv y)\Rightarrow x=y)$.
\end{teorema}

\begin{description}
\item[\sc Proof] For any $x$, $\mathbb{Z}(x)\Rightarrow \mathbb{M}_{\mathfrak{a}}^{\dagger}(x)$ (from Theorem \ref{todoZehMenge} and the fact that $\nu^{\mathfrak{a}}_{\dagger} \subset \nu^{\mathfrak{a}}$). But, from item {\sc i} of Definition \ref{definindoindiscerniveis}, for any ZF-set $x$ and any $y$, $x\equiv y \Leftrightarrow x=y$. Now, let $x$ be a Menge who is not a ZF-set, recalling $\mathfrak{a}$ satisfies all conditions imposed in Definition \ref{definindoindiscerniveis}. First we consider the case where $x$ acts only on atoms. Since no Menge $x$ in $\nu^{\mathfrak{a}}_{\dagger}$ acts on any atom $r$ without acting on $s$ when $r\equiv s$, then item {\sc iii} of Definition \ref{definindoindiscerniveis} grants $x\equiv y \Leftrightarrow x=y$. Finally, if $x$ and $y$ are any Mengen in $\nu^{\mathfrak{a}}_{\dagger}$, then $\mbox{TC}(x\big|_{\rank(t) = r \vee \rank(t) = s})\big|_{\mathfrak{a}[m]} \equiv \mbox{TC}(y\big|_{\rank(t) = r \vee \rank(t) = s})\big|_{\mathfrak{a}[m]}$ iff
$\mbox{TC}(x\big|_{\rank(t) = r \vee \rank(t) = s})\big|_{\mathfrak{a}[m]} = \mbox{TC}(y\big|_{\rank(t) = r \vee \rank(t) = s})\big|_{\mathfrak{a}[m]}$, where all occurrences of $m$ in formula $\mathfrak{a}[m]$ are free. So, item {\sc v} from Definition \ref{definindoindiscerniveis} grants $x\equiv y \Leftrightarrow x=y$.
\end{description}

\begin{teorema}\label{utilparalema7}
$\forall f(\mathbb{M}_{\mathfrak{a}}^{\dagger}(f) \Leftrightarrow (\sigma_f\neq\underline{\mathfrak{0}}\wedge f\subseteq\underline{\mathfrak{1}}\wedge \forall x(f[x]\Rightarrow (\forall t(t\equiv x\Rightarrow f[t])\wedge (\mathfrak{a}[x]\vee\mathbb{M}_{\mathfrak{a}}^{\dagger}(x))))))$.
\end{teorema}

\begin{description}
\item[\sc Proof] Analogous to the proof of Theorem \ref{teoremaquecapturantigadefinicaodezfconjunto} and a special case of Theorem \ref{tambemutilparalema7}.
\end{description}

Next theorem is a key result.

\begin{teorema}\label{nogo}
Any $f:x\to y$, such that $\mathbb{M}_{\mathfrak{a}}^{\dagger}(x)$ and $\mathbb{M}_{\mathfrak{a}}^{\dagger}(y)$, corresponds to a ZFU-function $f'$ where $\mathbb{M}_{\mathfrak{a}}^{\dagger}(f')$ iff neither $x$ nor $y$ acts on any atom.
\end{teorema}

\begin{description}
\item[\sc Proof] From Theorem \ref{funcaoZFfuncaocorrespondem} and Definition \ref{correspondencia} we know any $f:x\to y$, such that $\mathbb{M}_{\mathfrak{a}}^{\dagger}(x)$ and $\mathbb{M}_{\mathfrak{a}}^{\dagger}(y)$, corresponds to a ZFU-function where $Dom^{\mathfrak{F}}_f = Dom_{f'} = x$ and $Im^{\mathfrak{F}}_f = Im_{f'} = y$ (Theorem \ref{compartilhamdominioeimagem}) and $\mathbb{M}_{\mathfrak{a}}(f')$. The question is whether $\mathbb{M}_{\mathfrak{a}}^{\dagger}(f')$ or not. For the $\Leftarrow$ part, suppose $x$ acts on an atom $u$ from $\mathfrak{a}$. That entails $f(u) = v$ for some $v\neq\underline{\mathfrak{0}}$. But $f$ is in correspondence to a ZFU-function $f'$ iff $\forall u\forall v(f[u]\wedge f(u) = v)\Leftrightarrow (f'[(u,v)]\wedge f'((u,v)) = (u,v))$, where $x$, $y$, and $f'$ are Mengen in $\nu^{\mathfrak{a}}_{\dagger}$. Thus, any ordered pair $(u,v)$ is supposed to be a Menge in $\nu^{\mathfrak{a}}_{\dagger}$, where $f'$ acts. Nevertheless, according to Definition \ref{orderedpair}, $(u,v)$ is a function such that $\alpha\marcio{(u,v)}\alpha$ and $\beta\marcio{(u,v)}\beta$, where $u\marcio{\alpha} u$, and $u\marcio{\beta}u$ and $v\marcio{\beta} v$. Although $\alpha$ is a Menge in $\nu^{\mathfrak{a}}$, it is no Menge in $\nu^{\mathfrak{a}}_{\dagger}$. After all, any Menge from $\nu^{\mathfrak{a}}_{\dagger}$ who acts on any atom $u$ is supposed to act on $u'$ as well, where $u\equiv u'$. Since any atom $u$ in $\mathfrak{a}$ admits an $u'$ such that $u\neq u'$ and $u\equiv u'$, then $\neg\mathbb{M}_{\mathfrak{a}}^{\dagger}(\alpha)$. Therefore, $\neg\mathbb{M}_{\mathfrak{a}}^{\dagger}(f')$. An analogous rationale holds for the case where $y$ acts on any atom. For the $\Rightarrow$ part, suppose $f:x\to y$ corresponds to a ZFU-function $f'$ such that $\mathbb{M}_{\mathfrak{a}}^{\dagger}(f')$. That entails $f'$ acts on ordered pairs $(u,v)$ where $x[u]$ and $y[v]$. If either $u$ or $v$ is an atom, we get the same contradiction we got in the $\Leftarrow$ part.
\end{description}

\begin{teorema}\label{PPvale}
The Partition Principle holds in $\nu^{\mathfrak{a}}_{\dagger}$.
\end{teorema}

\begin{description}
\item[\sc Proof] Let $f'$ be a ZFU-function such that $\mathbb{M}_{\mathfrak{a}}^{\dagger}(f')$ and whose domain and range are, respectively, $x$ and $y$ (Definition \ref{dominioerangedezfu}). Then, there is a unique $f:x\to y$ in correspondence to $f'$ (Theorems \ref{funcaoZFfuncaocorrespondem} and \ref{compartilhamdominioeimagem}). Axiom {\sc F11} grants an injection $c:y\to x'$ such that $x'\subseteq x$. But since $f:x\to y$ is in correspondence to $f'$ where $\mathbb{M}_{\mathfrak{a}}^{\dagger}(f')$, that entails neither $x$ nor $y$ act on any atom (Theorem \ref{nogo}). Therefore, $c:y\to x'$ is in correspondence to a ZFU-function $c'$ where $\mathbb{M}_{\mathfrak{a}}^{\dagger}(c')$ (recall $x'\subseteq x$). Consequently, for any ZFU-function $f'$ in $\nu^{\mathfrak{a}}_{\dagger}$ such that $f$ is correspondence to $f'$, Partition Principle holds for $f'$ in $\nu^{\mathfrak{a}}_{\dagger}$, in the sense that $c'$ is a Menge in $\nu^{\mathfrak{a}}_{\dagger}$. Observe as well that $c'$ is injective in the usual sense concerning functions as collections of ordered pairs (labeled here as ZFU-functions).
\end{description}

\begin{teorema}\label{AEfalha}
Axiom of Choice fails in $\nu^{\mathfrak{a}}_{\dagger}$.
\end{teorema}

\begin{description}
\item[\sc Proof] Let $f:\mathfrak{a}\to r$ be a function where $r$ is an ordinal (Definition \ref{novadefinicaodeordinais}). Observe we have $\mathbb{M}_{\mathfrak{a}}^{\dagger}(\mathfrak{a})$ and $\mathbb{M}_{\mathfrak{a}}^{\dagger}(r)$. According to Theorem \ref{nogo}, there is no ZFU-function $f'$ such that $\mathbb{M}_{\mathfrak{a}}^{\dagger}(f')$ and $f$ is in correspondence to $f'$. Therefore, even if there is an injection $f$ from $\mathfrak{a}$ to $r$, there is no ZFU-function $f'$ in $\nu^{\mathfrak{a}}_{\dagger}$ which corresponds to such an injection. Therefore, $\mathfrak{a}$ cannot be well ordered in $\nu^{\mathfrak{a}}_{\dagger}$ by means of any ZFU-function $f'$.
\end{description}

The trick in the last two theorems is the fact that PP is a conditional: {\em if\/} there is a function from $x$ to $y$, then there is an injection from $y$ to $x'\subseteq x$. Nevertheless, certain ZFU-functions (those whose domains or ranges act on atoms) do not exist in $\nu^{\mathfrak{a}}_{\dagger}$ (Theorem \ref{nogo}).

Obviously, theorems \ref{PPvale} and \ref{AEfalha} do not grant we have succeeded into the introduction of a model of ZFU where PP holds but AC fails. After all, we still need to prove the remaining axioms of ZFU are, in some sense, true in our L\'evy universe. That is why we need the next Section. But before we go any further, we need to state our model, based on our previous discussions.

\begin{definicao}\label{pertencecomindicea}
$x\in_{\mathfrak{a}} y$ iff $y[x]\wedge \mathbb{M}_{\mathfrak{a}}^{\dagger}(y) \wedge (\exists t(t\equiv x\wedge (\mathfrak{a}[t]\vee\mathbb{M}_{\mathfrak{a}}^{\dagger}(t)))\Rightarrow y[t])$. If $\neg(x\in_{\mathfrak{a}} y)$, we abbreviate that as $x\not\in_{\mathfrak{a}} y$.
\end{definicao}

For example, formula $\alpha_r\in_{\mathfrak{a}} y$ is equivalent to $\overline{\alpha_r}\in_{\mathfrak{a}} y$ despite the fact that $\alpha_r \neq \overline{\alpha_r}$ for any finite ordinal $r$. Observe Theorem \ref{colecoesindiscerniveisaoiguais} grants we do not need to worry about any indiscernible of $x$ when $x\in_{\mathfrak{a}} y$ and $x$ is a Menge. In that case $x\in_{\mathfrak{a}} y$ is equivalent to $y[x]\wedge \mathbb{M}_{\mathfrak{a}}^{\dagger}(y)$. Thus,

\begin{teorema}\label{coisaesquisita}
Let $z$ be a function such that $\mathbb{M}_{\mathfrak{a}}^{\dagger}(z)$, and $x$ and $y$ be both atoms (or both Mengen) of $\nu^{\mathfrak{a}}_{\dagger}$. If $x\equiv y$, then $x\in_{\mathfrak{a}} z \Leftrightarrow y\in_{\mathfrak{a}} z$.
\end{teorema}

\begin{description}
\item[\sc Proof] Straightforward from Definition \ref{pertencecomindicea}.
\end{description}

If, e.g., $x$ is a Menge of $\nu^{\mathfrak{a}}_{\dagger}$, given by $\alpha_r\marcio{x}\alpha_r$ and $\overline{\alpha_r}\marcio{x}\overline{\alpha_r}$, then the only restrictions of $x$ who are Menge of $\nu^{\mathfrak{a}}_{\dagger}$ are $x$ itself and $\varphi_0$. In other words, atoms in a L\'evy universe always `come in pairs'. That is feasible thanks to the fact that the Axiom of Extensionality of ZFU cannot decide who is who between two indiscernible atoms. So, if we need any atom to form a Menge, we always take its `twin' with it. That is a different approach if we we compare it with permutation models \cite{Jech-03}. Finally,

\begin{definicao}\label{definindomodelodelevy}
Let $u = \underline{\mathfrak{1}}\big|_{\mathbb{M}_{\mathfrak{a}}^{\dagger}\,(t)}$, $F(t)$ be the formula $\exists a\exists b(a\in_{\mathfrak{a}} b)\Leftrightarrow t = (a,b)$, and $\epsilon = \underline{\mathfrak{1}}\big|_{F(t)}$. Then, $\mathfrak{M} = (u,\mathfrak{a},\varphi_0,\epsilon)$ is the {\em L\'evy model\/}.
\end{definicao}

Concerning last definition, it is worth to remark that $\epsilon$ acts on some terms where $u$ does not, namely, ordered pairs $(a,b)$ where $a$ is an atom and $b$ is a Menge, among others (see Theorem \ref{nogo}). In the next Section (specifically in Subsection \ref{venn}), it will be clear the meaning of $\exists x(P)$, where $P$ is a formula from ZFU's language. It means either $x$ is an atom ($\mathfrak{a}[x]$) or $x$ is a Menge in $\nu^{\mathfrak{a}}_{\dagger}$ ($\mathbb{M}_{\mathfrak{a}}^{\dagger}(x)$, which is equivalent to $u[x]$).

\section{ZFU is immersed in Flow}\label{standard}

We prove here vast portions of standard mathematics can be replicated within Flow.

\subsection{ZFU}

ZFU is a first-order theory with identity and one binary predicate letter $A_1^2$, such that the formula $A_1^2(x,y)$ is abbreviated as $x\in y$, if $x$ and $y$ are terms, and is read as `$x$ belongs to $y$' or `$x$ is an element of $y$'. The negation $\neg(x\in y)$ is abbreviated as $x\not\in y$. Besides, there are two constants, namely, $\emptyset$ and $A$, as primitive concepts. From this we are able to define a predicate $S$ as it follows: $S(a)$ iff $a\not\in A$. If $a\in A$ we say $a$ is an {\em atom\/}. Otherwise, we say $a$ is a {\em set\/}. The postulates of ZFU are the following:

\begin{description}
\item[\sc ZFU1 - Atoms] $\forall z(z\in A\Leftrightarrow (z\neq\emptyset \wedge \nexists x(x\in z)))$.

\item[\sc ZFU2 - Extensionality] $\forall_S x\forall_S y(\forall z(z\in x \Leftrightarrow z\in y)\Rightarrow x = y)$.

\item[\sc ZFU3 - Empty set] $\forall y(\neg(y\in \emptyset))$.

\item[\sc ZFU4 - Pair] $\forall x\forall y\exists_S z\forall t(t\in z\Leftrightarrow t = x \vee t = y)$.

\end{description}

$x\subseteq y =_{def} S(x)\wedge S(y)\wedge\forall z(z\in x\Rightarrow z\in y)$.

\begin{description}
\item[\sc ZFU5 - Power set] $\forall_S x\exists_S y\forall_S z(z\in y\Leftrightarrow z\subseteq x)$.
\end{description}

If $F(x)$ is a formula where there is no free occurrences of $y$, then:

\begin{description}
\item[\sc ZFU6$_F$ - Separation] $\forall_S z\exists_S y\forall x(x\in y \Leftrightarrow x\in z\wedge F(x))$.
\end{description}

The set $y$ is denoted by $\{x\in z\mid F(x)\}$.

If $\alpha(x,y)$ is a formula where all occurrences of $x$ and $y$ are free, then:

\begin{description}

\item[\sc ZFU7$_\alpha$ - Replacement] $\forall x\exists!y\alpha(x,y)\Rightarrow\forall_S z\exists_S w\forall t(t\in w
\Leftrightarrow\exists s(s\in z\wedge\alpha(s,t)))$.

\item[\sc ZFU8 - Union set] $\forall_S x\exists_S y\forall z (z\in y\Leftrightarrow\exists t(z\in t\wedge t\in x))$.

\end{description}

The set $y$ from {\sc ZF7} is abbreviated as $y = \bigcup_{t\in x} t$. Intersection $\bigcap$ is defined from union $\bigcup$ in the usual way. We adopt the convention that $x\cap y = \emptyset$ if either $x$ or $y$ is an atom.

\begin{description}
\item[\sc ZFU9 - Infinity] $\exists_S x(\emptyset\in x\wedge \forall_S y(y\in x\Rightarrow y\cup\{y\}\in x))$.
\end{description}

\begin{description}
\item[\sc ZFU10 - Choice] $\forall_S x(\forall_S y\forall_S z((y\in x\wedge z\in x\wedge y\neq z)\Rightarrow (y\neq\emptyset \wedge
y\cap z = \emptyset))\Rightarrow \exists_S y\forall_S z(z\in x\Rightarrow \exists w (y\cap z = \{w\})))$.

\item[\sc ZFU11 - Regularity] $\forall_S x(x\neq\emptyset\Rightarrow \exists y(y\in x \wedge x\cap y = \emptyset))$.

\end{description}

\subsection{L\'evy model}\label{venn}

\begin{definicao}\label{tabelaqueinteressa}
Consider the next translation table,

\begin{center}
\begin{tabular}{|c|c|}\hline
\multicolumn{2}{|c|}{\sc Translating ZFU into $\mathfrak{M}$}\\ \hline \multicolumn{1}{|c|}{\sc ZFU} & \multicolumn{1}{|c|}{$\mathfrak{M}$}\\ \hline\hline
$\forall x$ & $\forall_{\mathbb{A}_{\mathfrak{a}}\;\vee \mathbb{M}_{\mathfrak{a}}^{\dagger}}\,x$\\
$\exists x$ & $\exists_{\mathbb{A}_{\mathfrak{a}}\;\vee \mathbb{M}_{\mathfrak{a}}^{\dagger}}\,x$\\
$S(x)$ & $\mathbb{M}_{\mathfrak{a}}^{\dagger}(x)$\\
$A$ & $\mathfrak{a}$\\
$\emptyset$ & $\varphi_0$\\
$x\in y$ & $x\in_{\mathfrak{a}}y$\\
$x\subseteq y$ & $\forall t(t\in_{\mathfrak{a}}x\Rightarrow t\in_{\mathfrak{a}}y)$\\
$x = y$ & $x = y$\\ \hline
\end{tabular}
\end{center}
\noindent
where $\mathfrak{M}$ is the L\'evy model and $\mathbb{A}_{\mathfrak{a}}$ is a monadic predicate where $\mathbb{A}_{\mathfrak{a}}(x)$ says $x$ is an atom ($\mathbb{A}(x)$) and $\mathfrak{a}[x]$. If $\Xi$ is a formula from ZFU, then its translation by means of the table above is denoted by $Translated(\Xi)$. Besides,

$$\mathfrak{M}\vDash \Xi\;\;\;\mbox{iff}\;\;\;\vdash_{\mbox{\boldmath{$\mathfrak{F}$}}}Translated(\Xi),$$
\noindent
where $\mbox{\boldmath{$\mathfrak{F}$}}$ is Flow theory. We read $\mathfrak{M}\vDash \Xi$ as `$\Xi$ {\em is true in\/} $\mathfrak{M}$'.
\end{definicao}

\begin{proposicao}\label{preservamatematica}
All axioms of ZFU are true in the L\'evy model, with the only exception of the Axiom of Choice.
\end{proposicao}

The proof of Proposition \ref{preservamatematica} is made through the following lemmas. For the sake of readability, translated quantifiers $\forall_{\mathbb{A}_{\mathfrak{a}}\;\vee \mathbb{M}_{\mathfrak{a}}^{\dagger}}\,x$ and $\exists_{\mathbb{A}_{\mathfrak{a}}\;\vee \mathbb{M}_{\mathfrak{a}}^{\dagger}}\,x$ will be simply written as $\forall x$ and $\exists x$, respectively. There is no risk of confusion since the proofs of all lemmas refer to translated formulas. Observe as well that $\forall_S x$ and $\exists_S x$ from ZFU language are translated as $\forall_{\mathbb{M}_{\mathfrak{a}}^{\dagger}}\,x$ and $\exists_{\mathbb{M}_{\mathfrak{a}}^{\dagger}}\,x$, respectively.

\begin{lema}
$\mathfrak{M}\vDash \mbox{\sc ZFU1}$.
\end{lema}

\begin{description}
\item[\sc Proof] $Translated(\mbox{ZFU1})$ is $\forall z(z\in_{\mathfrak{a}}\mathfrak{a} \Leftrightarrow (z\neq\varphi_0\wedge \nexists x(x\in_{\mathfrak{a}} z)))$. For the $\Rightarrow$ part, observe $z\in_{\mathfrak{a}}\mathfrak{a}$ is equivalent to say $z$ is an atom (Definitions \ref{definindoatomos} and \ref{tijoloapropriado}). Since no atom is a restriction of $\underline{\mathfrak{1}}$, then any atom $z$ is different of $\varphi_0$. Besides, no term belongs to (we are talking about $\in_{\mathfrak{a}}$) any atom $z$, since no atom is a Menge. For the $\Leftarrow$ part, observe $z$ is either an atom or a Menge. We already know there is no $x$ such that $x\in_{\mathfrak{a}}\varphi_0$, besides the fact that $\mathbb{M}_{\mathfrak{a}}^{\dagger}(\varphi_0)$. Suppose $z\neq\varphi_0$ is an empty Menge in the L\'evy universe, i.e., $\nexists x(x\in_{\mathfrak{a}} z)$. That entails $z$ does not act on any Menge nor on any atom of the L\'evy universe. If $z$ acts on any other term, then $z$ is no Menge in $\nu^{\mathfrak{a}}_{\dagger}$. Thus, $\varphi_0$ is the only empty Menge in $\nu^{\mathfrak{a}}_{\dagger}$. Therefore, $z$ is an atom.
\end{description}

\begin{lema}\label{lemaprimeiro}
$\mathfrak{M}\vDash \mbox{\sc ZFU2}$.
\end{lema}

\begin{description}
\item[\sc Proof] $Translated(\mbox{\sc ZFU2})$ is $\forall_{\mathbb{M}_{\mathfrak{a}}^{\dagger}}\,x\forall_{\mathbb{M}_{\mathfrak{a}}^{\dagger}}\,y(\forall z(z\in_{\mathfrak{a}} x \Leftrightarrow z\in_{\mathfrak{a}} y)\Rightarrow x = y)$. If $x$ and $y$ are Mengen and $z\in_{\mathfrak{a}} x$ and $z\in_{\mathfrak{a}} y$, then $x(z) = z$ and $y(z) = z$. Besides, if $z'\equiv z$, then $x(z') = z'$ and $y(z') = z'$ (Theorem \ref{coisaesquisita}). If either $z\not\in_{\mathfrak{a}} x$ or $z\not\in_{\mathfrak{a}} x$, then either $x(z) = \underline{\mathfrak{0}}$ or $y(z) = \underline{\mathfrak{0}}$ or $z = x$ or $z = y$ (where we may replace all occurrences of $z$ by $z'$ if $z\equiv z'$). So, $Translated(\mbox{ZFU2})$ considers the case where both $x$ and $y$ share all their images $x(z)$ and $y(z)$ for any $z$, except perhaps for $z = x$ or $z = y$. That means $x\sim y$ (Definition \ref{definesim}). But Theorem \ref{marcionoslibertou} says no emergent function has any clone which is emergent. Thus, $Translated(\mbox{ZFU2})$ is a theorem of Flow thanks to Theorem  \ref{igualdadefuncoes}.
\end{description}

\begin{lema}
$\mathfrak{M}\vDash \mbox{\sc ZFU3}$.
\end{lema}

\begin{description}
\item[\sc Proof] $Translated(\mbox{\sc ZFU3})$ is $\forall y(\neg(y\in_{\mathfrak{a}} \varphi_0))$. But $\mathbb{M}_{\mathfrak{a}}^{\dagger}(\varphi_0)$, which means $\varphi_0$ is a Menge in the L\'evy universe. Since $\varphi_0$ does not act on any term, then there is no $y$ such that $y\in_{\mathfrak{a}} \varphi_0$.
\end{description}

\begin{lema}\label{lemadoparordenado}
$\mathfrak{M}\vDash \mbox{\sc ZFU4}$.
\end{lema}

\begin{description}
\item[\sc Proof] $Translated(\mbox{\sc ZFU4})$ is $\forall x\forall y\exists_{\mathbb{M}_{\mathfrak{a}}^{\dagger}}\, z\forall t(t\in_{\mathfrak{a}} z\Leftrightarrow t = x \vee t = y)$. Let $h = \varphi_4\big|^{\alpha(t,r)}$, where $\alpha(t,r)$ is the next formula: $(t = \varphi_0\Rightarrow r = x)\wedge ((t = \varphi_1 \wedge \exists x'(x'\neq x\wedge x'\equiv x))\Rightarrow r = x')\wedge ((t = \varphi_1\wedge \nexists x'(x'\neq x\wedge x'\equiv x))\Rightarrow r = \underline{\mathfrak{0}})\wedge (t = \varphi_2\Rightarrow r = y)\wedge ((t = \varphi_3\wedge \exists y'(y'\neq y\wedge y'\equiv y))\Rightarrow r = y')\wedge ((t = \varphi_3\wedge \nexists y'(y'\neq y\wedge y'\equiv y))\Rightarrow r = \underline{\mathfrak{0}})\wedge ((t\neq\varphi_0 \wedge t\neq\varphi_1 \wedge t\neq\varphi_2 \wedge t\neq \varphi_3)\Rightarrow r = \underline{\mathfrak{0}})$. Observe the quantifier $\exists$ in formula $\alpha(t,r)$ is supposed to be read as in the table from Definition \ref{venn}. According to {\sc F10}$_{\alpha}$, $\sigma_h\neq\underline{\mathfrak{0}}$ and $\sigma(Im_h^{\mathfrak{F}})\neq\underline{\mathfrak{0}}$. Now, let $z = Im_h^{\mathfrak{F}} = \underline{\mathfrak{1}}\big|_{t = h(\varphi_0)\vee t = h(\varphi_1)\vee t = h(\varphi_2)\vee t = h(\varphi_3)}$. From Definition \ref{levy} and Theorem \ref{utilparalema7}, $\mathbb{M}_{\mathfrak{a}}^{\dagger}(z)$. Thus, $t\in_{\mathfrak{a}} z\Leftrightarrow t = x \vee t = y$.
\end{description}

If the reader wants a more intuitive view about this last lemma, observe the only terms $x$ and $y$ in a Levy universe who are indiscernible and different from each other are atoms, according to Theorem \ref{colecoesindiscerniveisaoiguais}. Thus, any pair $z$ of Mengen has either one or two elements with respect to $\in_{\mathfrak{a}}$, as it is expected. Now, let $x$ be an atom. Hence, the pair $z$ such that $t\in_{\mathfrak{a}} z\Leftrightarrow t = x$ is a Menge who acts on $x$ and its conjugate $x'$. We say $x'$ is a conjugate of $x$ iff $x\equiv x'$ and $x\neq x'$ (see Theorem \ref{coisaesquisita}). As a final example, if $x = \alpha_{\varphi_0}$ and $y = \overline{\alpha_{\varphi_1}}$, then the pair $z$ acts on four terms, namely, $\alpha_{\varphi_0}$, $\overline{\alpha_{\varphi_0}}$, $\alpha_{\varphi_1}$, and $\overline{\alpha_{\varphi_1}}$.

\begin{lema}
$\mathfrak{M}\vDash \mbox{\sc ZFU5}$.
\end{lema}

\begin{description}
\item[\sc Proof] $Translated(\mbox{\sc ZFU5})$ is $\forall_{\mathbb{M}_{\mathfrak{a}}^{\dagger}}\, x\exists_{\mathbb{M}_{\mathfrak{a}}^{\dagger}}\, y\forall_{\mathbb{M}_{\mathfrak{a}}^{\dagger}}\, t(t\in_{\mathfrak{a}} y\Leftrightarrow t\subseteq x)$. Let $y = \underline{\mathfrak{1}}\big|_{t\subseteq x}$ (see the translation of $\subseteq$ in Definition \ref{tabelaqueinteressa}). According to {\sc F8}, every term where $y$ acts is emergent. Besides, $y$ itself is emergent (which entails $\sigma_y\neq\underline{\mathfrak{0}}$). More than that, since $\mathbb{M}_{\mathfrak{a}}^{\dagger}(x)$, then $\mathbb{M}_{\mathfrak{a}}^{\dagger}(y)$ (Definitions \ref{levy} and \ref{tijolo} grant every restriction of a Menge is a Menge). Recall Definition \ref{tijolo} makes use exactly of the notion of restricted power $\wp$. Finally, observe $t\in_{\mathfrak{a}} y\Leftrightarrow t\subseteq x$, in the sense that $y[t]\wedge \mathbb{M}_{\mathfrak{a}}^{\dagger}(y)$ is equivalent to $t\in_{\mathfrak{a}} y$ according to Definition \ref{pertencecomindicea} and its subsequent paragraph.
\end{description}

\begin{lema}
For any formula $F$ with the syntactic conditions in {\sc ZFU6}$_F$ we have $\mathfrak{M}\vDash \mbox{\sc ZFU6}_F$.
\end{lema}

\begin{description}
\item[\sc Proof] $Translated(\mbox{\sc ZFU6}_F)$ is $\forall_{\mathbb{M}_{\mathfrak{a}}^{\dagger}}\, z\exists_{\mathbb{M}_{\mathfrak{a}}^{\dagger}}\, y\forall t(t\in_{\mathfrak{a}} y \Leftrightarrow t\in_{\mathfrak{a}} z\wedge F(t))$, if $F(t)$ is a formula with no free occurrences of $y$. Let $y = z\big|_{\exists r\exists s((r\equiv s\wedge F(s))\Leftrightarrow t = r)}$ (remember every term is indiscernible from itself and, eventually, it can be indiscernible from something else). Since $\mathbb{M}_{\mathfrak{a}}^{\dagger}(z)$ and $y\subseteq z$, then $\mathbb{M}_{\mathfrak{a}}^{\dagger}(y)$ (recall that if $y$ acts on any $t$, then it acts on any $t'$ such that $t\equiv t'$). Observe $t\in_{\mathfrak{a}} y$ is equivalent to $t\in_{\mathfrak{a}} z\wedge F(t)$, since Theorem \ref{coisaesquisita} grants $t\in_{\mathfrak{a}} y$ is equivalent to $s\in_{\mathfrak{a}} y$ when $t\equiv s$. Recalling Observation \ref{consideracoessobrerestricao}, the calculation of $y = z\big|_{\exists r\exists s((r\equiv s\wedge F(s))\Leftrightarrow t = r)}$ does  not take into account any $t = y$. So, any occurrence of $y$ in $F$ turns out to be irrelevant.
\end{description}

From an intuitive point of view, the only atypical situation in $t\in_{\mathfrak{a}} z\wedge F(t)$ takes place when $t$ is an atom. In that case, formula $t\in_{\mathfrak{a}} z\wedge F(t)$ can be read as `$t$ and its conjugate belong to $z$ and, at least for one of them, formula $F$ holds'.

\begin{lema}
For any formula $\alpha$ with the syntactic conditions in {\sc ZFU7}$_{\alpha}$ we have $\mathfrak{M}\vDash \mbox{\sc ZFU7}_{\alpha}$.
\end{lema}

\begin{description}
\item[\sc Proof] $Translated(\mbox{\sc ZFU7}_{\alpha})$ is $\forall x\exists!y\alpha(x,y)\Rightarrow\forall_{\mathbb{M}_{\mathfrak{a}}^{\dagger}}\, z\exists_{\mathbb{M}_{\mathfrak{a}}^{\dagger}}\, w\forall t(t\in_{\mathfrak{a}} w
\Leftrightarrow\exists s(s\in_{\mathfrak{a}} z\wedge\alpha(s,t)))$, where all occurrences of $x$ and $y$ in $\alpha(x,y)$ are free. Due to the use of quantifier $\exists!$, we are supposed to rewrite {\sc ZFU7} by means of $\forall$ and $\exists$ (as defined at the end of the second paragraph of Section \ref{secao2}), and only then translate it. Let $g = z\big|^{\alpha(s,t)}$ be defined from {\sc F10}$_{\alpha}$ by formula $\alpha(s,t)$ provided by $Translated(\mbox{\sc ZFU7}_{\alpha})$. From {\sc F10}$_{\alpha}$, $\sigma_g\neq\underline{\mathfrak{0}}$ and $\sigma(Im_g^{\mathfrak{F}})\neq\underline{\mathfrak{0}}$. That makes $g$ emergent. Now, let $w' = Im_g^{\mathfrak{F}} = \underline{\mathfrak{1}}\big|_{\exists s(z[s]\wedge t = g(s))}$. But $w'$ is a restriction of $\underline{\mathfrak{1}}$, with $\sigma_{w'}\neq\underline{\mathfrak{0}}$, who acts only on either atoms or Mengen. Therefore, from Theorem \ref{tambemutilparalema7}, $w'$ is a Menge in the appropriate atomic von Neumann universe $\nu^{\mathfrak{a}}$. Now, let $w = \underline{\mathfrak{1}}\big|_{\exists s(w'[s]\Leftrightarrow (t\equiv s\wedge (\mathfrak{a}[t]\vee\mathbb{M}_{\mathfrak{a}}(t))))}$. Observe $w'$ is a restriction of $w$. Since $\mathbb{M}_{\mathfrak{a}}(w')$, so is $w$, where, by the way, both share the same rank (Definitions \ref{tijolo} and \ref{rankeandocomatomos}). But, more than that, $\mathbb{M}_{\mathfrak{a}}^{\dagger}(w)$ (Definition \ref{levy}), since for any $s$ where $w$ acts, the latter acts on any indiscernible of $s$ as well. Thus, for any $t$ (either an atom or a Menge in $\nu^{\mathfrak{a}}_{\dagger}$), $t\in_{\mathfrak{a}} w
\Leftrightarrow\exists s(s\in_{\mathfrak{a}} z\wedge\alpha(s,t))$.
\end{description}

Nothing prevents formula $\alpha(s,t)$ from last lemma be defined in a way such that $t$ is an atom (for a given $s$) but there is no $s'$ where $\alpha(s',t')$ and $t'\equiv t$ and $t'\neq t$. That is the reason why we had to use $\nu^{\mathfrak{a}}$ before we get to $\nu^{\mathfrak{a}}_{\dagger}$.

\begin{lema}\label{lemadauniaoZF}
$\mathfrak{M}\vDash \mbox{\sc ZFU8}$.
\end{lema}

\begin{description}
\item[\sc Proof] $Translated(\mbox{\sc ZFU8})$ is $\forall_{\mathbb{M}_{\mathfrak{a}}^{\dagger}}\, x\exists_{\mathbb{M}_{\mathfrak{a}}^{\dagger}}\, y\forall t (t\in_{\mathfrak{a}} y\Leftrightarrow\exists z(t\in_{\mathfrak{a}} z\wedge z\in_{\mathfrak{a}} x))$. Let $y = \bigcup_{z\in_{\mathfrak{a}}\; x} z =_{def} \underline{\mathfrak{1}}\big|_{\exists z(z\in_{\mathfrak{a}}\; x\wedge t\in_{\mathfrak{a}}\; z)}$. Axiom {\sc F8} entails $\sigma_y\neq\underline{\mathfrak{0}}$. Theorem \ref{utilparalema7} entails $\mathbb{M}_{\mathfrak{a}}^{\dagger}(y)$.
\end{description}

\begin{lema}
$\mathfrak{M}\vDash \mbox{\sc ZFU9}$.
\end{lema}

\begin{description}
\item[\sc Proof] $Translated(\mbox{\sc ZFU9})$ is $\exists_{\mathbb{M}_{\mathfrak{a}}^{\dagger}}\, x(\varphi_0\in_{\mathfrak{a}} x\wedge \forall_{\mathbb{M}_{\mathfrak{a}}^{\dagger}}\, y(y\in_{\mathfrak{a}} x\Rightarrow \sigma_y\in_{\mathfrak{a}} x))$. Let $x = \omega$, from Definition \ref{defineordinalimite}. Recall any ZF-set is a Menge in the L\'evy universe (Theorem \ref{todoZehMenge}).
\end{description}

\begin{lema}
$\mathfrak{M}\vDash \mbox{\sc ZFU11}$.
\end{lema}

\begin{description}
\item[\sc Proof] $Translated(\mbox{\sc ZFU11})$ is $\forall_{\mathbb{M}_{\mathfrak{a}}^{\dagger}}\, x(x\neq\varphi_0\Rightarrow \exists y(y\in_{\mathfrak{a}} x \wedge x\cap y = \varphi_0))$. Any Menge $x$ in the L\'evy universe has a rank (Definition \ref{rankeandocomatomos}). Let $y\in_{\mathfrak{a}} x$ be any term such that there is no $z$ where $\rank(y)[\rank(z)]$ and $z\in_{\mathfrak{a}} x$ (see Theorem \ref{teoremasobreacoeserank}, which can be extended to any atomic von Neumann universe), i.e., $y$ has the least rank among the elements of $x$ (relative to $\in_{\mathfrak{a}}$). Then, there is not $t$ such that $t\in_{\mathfrak{a}} x$ and $t\in_{\mathfrak{a}} y$.
\end{description}

This concludes the proof of Proposition \ref{preservamatematica}. Hence, Proposition \ref{preservamatematica} and Theorems \ref{PPvale} and \ref{AEfalha} guarantee  $\mathfrak{M}$ is a model of ZFU where the Partition Principle holds but the Axiom of Choice fails. After all, $\mathfrak{M}$ (Definition \ref{definindomodelodelevy}) is defined from $\nu^{\mathfrak{a}}_{\dagger}$.

\subsection{Other models}

\begin{definicao}
Let $u = \underline{\mathfrak{1}}\big|_{\mathbb{Z}(t)}$, $F(t)$ be the formula $\exists a\exists b(a\in_{\mathbb{Z}} b)\Leftrightarrow t = (a,b)$, and $\epsilon = \underline{\mathfrak{1}}\big|_{F(t)}$. Then, $\mathfrak{N} = (u,\epsilon)$ is the {\em von Neumann model\/}, where $\mathbb{Z}$ is the predicate ``to be a ZF-set'' from Definition \ref{defineZFconjunto} and $\in_{\mathbb{Z}}$ is given in Definition \ref{pertencercomoagir}.
\end{definicao}

\begin{definicao}\label{outrosmodelos}
Consider the next translation table.
\begin{center}
\begin{tabular}{|c|c|}\hline
\multicolumn{2}{|c|}{\sc Translating ZF into $\mbox{$\mathfrak{N}$}$}\\ \hline \multicolumn{1}{|c|}{\sc ZF} & \multicolumn{1}{|c|}{$\mbox{$\mathfrak{N}$}$}\\ \hline\hline
$\forall$ & $\forall_{\mathbb{Z}}$\\
$\exists$ & $\exists_{\mathbb{Z}}$\\
$x\in y$ & $x\in_{\mathbb{Z}} y$\\
$x\subseteq y$ & $x\subseteq y$\\ \hline
\end{tabular}
\end{center}
If $\Xi$ is a formula from ZF, then its translation by means of the table above is denoted by $Translated(\Xi)$. Besides,

$$\mathfrak{N}\vDash \Xi\;\;\;\mbox{iff}\;\;\;\vdash_{\mbox{\boldmath{$\mathfrak{F}$}}}Translated(\Xi),$$
\noindent
where $\mbox{\boldmath{$\mathfrak{F}$}}$ is Flow theory. We read $\mathfrak{N}\vDash \Xi$ as `$\Xi$ {\em is true in\/} $\mathfrak{N}$'.
\end{definicao}

\begin{proposicao}
All axioms of ZF are true in the von Neumann model.
\end{proposicao}

\begin{description}
\item[\sc Proof] Analogous to the proof of Proposition \ref{preservamatematica}, but much simpler.
\end{description}

If the reader is still suspicious about last proof, observe $\nu\subset \nu^{\mathfrak{a}}_{\dagger} \subset \nu^{\mathfrak{a}}$.

To obtain a model of ZFC, all we have to do is to replace {\sc F11} by $Translated$(AC) according to Definition \ref{outrosmodelos}. Obviously, both AC and PP are true in such a model.

If we replace $\in_{\mathfrak{a}}$ by $\in$ (Definition \ref{pertinenciasemdecisao}) and all occurrences of $\mathbb{M}_{\mathfrak{a}}^{\dagger}$ by $\mathbb{M}_{\mathfrak{a}}$ in Definition \ref{tabelaqueinteressa} (keeping {\sc F11}), we have a model of ZFU where we are unable to answer if AC is independent from PP.

\section{On the consistency of Flow}

It is far from us to answer if Flow is consistent, due to G\"odel's Second Incompleteness Theorem \cite{Godel-92}. But we can talk about the consistency of $\mbox{\boldmath{$\mathfrak{F}$}}$ relative to ZF, ZFC, and ZFU.

In last Section we proved that, if Flow is consistent, then ZF, ZFC, and ZFU are consistent. We did that by means of a metatheory which provided a translation from the language of ZF (and some of its variations) into Flow's language. Thus, there seems to be two possibilities, if we want to assess the consistency of Flow relative to ZF:

\begin{enumerate}

\item We introduce a translation from the language of Flow into the language of ZF and, then, prove every translated axiom of Flow is a theorem in ZF.

\item We use some other metatheory (quite different from a simple translation) to decide if Flow has a greater (or equal) consistency strength than ZF.

\end{enumerate}

The first possibility seems to bump into many obstacles. How to interpret $\underline{\mathfrak{0}}$ and $\underline{\mathfrak{1}}$ into the language of ZF? If $\underline{\mathfrak{0}}$ does not translate as the empty set, what is it? If $\underline{\mathfrak{0}}$ is to be interpreted as the empty set, who is $\varphi_0$? Recall $\underline{\mathfrak{0}}$ and $\varphi_0$ are the only clones in the present formulation of Flow. Hence, to avoid inconsistency with extensionality in ZF, we cannot interpret both $\underline{\mathfrak{0}}$ and $\varphi_0$ as the empty set. On the other hand, if $\underline{\mathfrak{1}}$ acts on every single term (except itself and $\underline{\mathfrak{0}}$), how to interpret it within ZF? It is worth to remark that $\underline{\mathfrak{1}}$ is not a proper class, in the sense of a theory like, e.g., NBG \cite{Mendelson-97}. After all, $\underline{\mathfrak{1}}$ acts on many universes (like $\nu$, $\nu^a$ for any brick of atoms $a$, $\nu^{\mathfrak{a}}_{\dagger}$, and so on) where those universes behave like proper classes. Thus, if we associate actions to membership relations in the sense of Definition \ref{pertencecomindicea}, that would entail many classes belonging to $\underline{\mathfrak{1}}$. Besides, our main results show the usual correspondence between sets and identity functions is not that trivial. What Theorem \ref{nogo} (and its consequences in Theorems \ref{PPvale} and \ref{AEfalha}) reveals is that a set $x$ is not necessarily associated to a function $f:x\to x$ (in the ordinary set-theoretic sense, where a function is a set of ordered pairs) such that $\forall t(f(t) = t)$. That non-correspondence takes place in our L\'evy universe when $x$ has atoms as elements.

Concerning the second possibility, could we try and use forcing \cite{Jech-03} over Flow? Obviously, no, at least for now. Forcing requires the previous existence of at least one model in order to conceive other models for the sake of investigation about relative consistency among formal theories. As remarked by Saharon Shelah \cite{Shelah-03}, ``forcing can be used only to make the universe `fatter', not `taller'. In technical terms: if we use forcing starting from models of ZFC to prove that `ZFC neither proves nor refutes statement A' (or equivalently, each of `ZFC + A' and `ZFC + non-A' is consistent), then the `consistency strengths' of `ZFC', `ZFC+A', `ZFC + non-A' are all equal''. The `universe' mentioned here refers to a previously known model. The `taller' independence results mentioned by Shelah refer to situations where the consistency strength of `ZFC + statement A' is strictly higher than the consistency strength of ZFC alone. That `statement A' could be, for example, `ZFC is consistent'. Such `statement A' is far beyond the reach of forcing techniques. But our concern here is a little less demanding. The point here is the fact that we have no metalinguistic interpretation which could serve as a model of Flow. That is a task to be carried on in the future. Besides, we admit the possibility that Flow demands some quite new techniques for investigating its consistency strength relatively to other theories like ZF.

\section{Final remarks}\label{finalremarks}

As a reference to Heraclitus's flux doctrine, we are inclined to refer to all terms of Flow as {\em fluents\/}, rather than functions. In this sense Flow can be roughly regarded as a {\em theory of fluents\/}. That is also an auspicious homage to the {\em Method of Fluxions\/} by Isaac Newton \cite{Newton-64}. The famous `natural philosopher' referred to functions as fluents, and their derivatives as {\em fluxions\/}. Whether Newton was inspired or not by Heraclitus, that is historically uncertain (\cite{Royal-46}, page 38). Notwithstanding, we find such a coincidence quite inspiring and utterly opportune.

From the mathematical point of view, our framework was strongly motivated by von Neumann's set theory \cite{vonNeumann-25}, as we briefly discussed in the Introduction. But what are the main differences between Flow and von Neumann's ideas?

We refer to von Neumann's set theory as $\mbox{\boldmath{$\mathfrak{N}$}}$. The first important difference lurks in the way how von Neumann seemed to understand functions: ``a function can be regarded as a set of pairs, and a set as a function that can take two values... the two notions are completely equivalent''. That philosophical viewpoint seems to be committed to a set-theoretic framework. Here we follow a quite different path, since our main result suggests those notions are not necessarily equivalent: {\sc i}) In $\mbox{\boldmath{$\mathfrak{N}$}}$ there are two privileged objects termed $A$ and $B$ which resemble our terms $\underline{\mathfrak{0}}$ and $\underline{\mathfrak{1}}$. Nevertheless, in $\mbox{\boldmath{$\mathfrak{N}$}}$ there is no further information about $A$ and $B$. Besides, those constants play a different role if we compare to $\underline{\mathfrak{0}}$ and $\underline{\mathfrak{1}}$ in $\mbox{\boldmath{$\mathfrak{F}$}}$, as we can see in the next item. {\sc ii}) In $\mbox{\boldmath{$\mathfrak{N}$}}$ there are two sorts of objects, namely, arguments and functions. Eventually some of those objects are both of them. And when an object $f$ is an argument and a function, which takes only values $A$ and $B$, then $f$ is a set. Our terms $\underline{\mathfrak{0}}$ and $\underline{\mathfrak{1}}$ have no similar role. Besides, we do not need to distinguish functions from any other kind of term. All objects of $\mbox{\boldmath{$\mathfrak{F}$}}$ are functions. {\sc iii}) Von Neumann believed a distinction between arguments and functions was necessary to avoid the well known antinomies from naive set theory. Nevertheless, we proved a simple axiom of Self-Reference ({\sc F2}) is enough to avoid such a problem. {\sc iv}) The distinction between sets and other collections which are `too big' to be sets depends on considerations if a specific term in $\mbox{\boldmath{$\mathfrak{N}$}}$ is both an argument and a function. Within $\mbox{\boldmath{$\mathfrak{F}$}}$ that distinction depends on considerations regarding $\mathfrak{F}$-successor. {\sc v}) Axiom I4 of $\mbox{\boldmath{$\mathfrak{N}$}}$ (\cite{vonNeumann-25}, page 399) says any function can be identified by its images. Our weak extensionality {\sc F1} allows us to derive a similar result as a non-trivial theorem (Theorem \ref{igualdadefuncoes}). {\sc vi}) One primitive concept in $\mbox{\boldmath{$\mathfrak{N}$}}$ is a binary functional letter (using modern terminology) which allows to define ordered pairs. In $\mbox{\boldmath{$\mathfrak{F}$}}$ that assumption is unnecessary. {\sc vii}) In $\mbox{\boldmath{$\mathfrak{N}$}}$ there are many constant functions (\cite{vonNeumann-25}, page 399, axiom II2). In $\mbox{\boldmath{$\mathfrak{F}$}}$ there are only two constant functions, namely, $\underline{\mathfrak{0}}$ and $\varphi_0$. {\sc viii}) In $\mbox{\boldmath{$\mathfrak{N}$}}$ the well-ordering theorem and the axiom of choice are consequences of its postulates. In $\mbox{\boldmath{$\mathfrak{F}$}}$ such a phenomenon does not take place. Actually, it makes no sense at all to state something like AC within Flow.

Some open problems refer to how precisely Flow is related to Lambda Calculus \cite{Barendregt-13}, Category Theory \cite{Hatcher-68} \cite{Lawvere-03}, String Diagrams \cite{Dixon-13}, and Autocategories \cite{Guitart-14}. Our technique for building a model of ZFU is supposed to be compared with permutation models as well \cite{Jech-03}. But, for now, the most important open problem is whether AC can be proven to be independent from PP within ZF.

\section{Acknowledgements}

We thank Aline Zanardini, Bruno Victor, and Cl\'eber Barreto for insightful discussions which inspired this work. We acknowledge with thanks as well Asaf Karagila, Hanul Jeon, Samuel Gomes da Silva, Newton da Costa, Jean-Pierre Marquis, Ed\'elcio Gon\c calves de Souza, D\'ecio Krause, Jonas Arenhart, Kherian Gracher, Bryan Leal Andrade, and B\'arbara Guerreira for valuable remarks and criticisms concerning earlier and problematic versions of this paper. AK was the most patient, critic, and thoughtful reader of a very poor earlier version of this work. We were also benefitted by comments and criticisms made by several participants of seminars delivered at Federal University of Paran\'a and Federal University of Santa Catarina. Finally we acknowledge Pedro D. Dam\'azio's provocation, more than 20 years ago, which worked as the first motivation towards the development of Flow.

\end{document}